\documentclass[10pt]{amsart}
\usepackage{a4}
\usepackage{amssymb}
\usepackage{array}
\usepackage{booktabs}
\usepackage{multirow}
\usepackage{hhline}

\newtheorem{thm}{Theorem}

\newtheorem{cor}[thm]{Corollary}

\newtheorem{rem}[thm]{Remark}

\numberwithin{thm}{section}
\numberwithin{equation}{section}
\newcommand{\Real}{\mathbb R}

\newcommand{\norm}[1]{\left\Vert#1\right\Vert}
\newcommand{\abs}[1]{\left\vert#1\right\vert}

\newcommand{\eps}{\varepsilon}

\newcommand{\A}{\mathcal{A}}

\newcommand{\C}{\mathcal{C}}

\newcommand{\E}{\mathbb{E}}

\newcommand{\vep}{\varepsilon}
\newcommand{\F}{\mathcal{F}}
\newcommand{\G}{\mathcal{G}}
\newcommand{\Hi}{\mathcal{H}}
\newcommand{\I}{\mathcal{I}}
\newcommand{\J}{\mathcal{J}}
\newcommand{\K}{\mathcal{K}}

\newcommand{\M}{\mathcal{M}}
\newcommand{\n}{\mathbb{N}}
\newcommand{\N}{\mathcal{N}}
\newcommand{\R}{\mathcal{R}}
\newcommand{\rb}{\textbf{r}}
\newcommand{\s}{\mathcal{S}}

\newcommand{\U}{\mathcal{U}}
\newcommand{\und}{\underline}

\newcommand{\Om}{\Omega}

\begin{document}
\title[$q$-Chaos]
{$q$-Chaos}
\author{Marius Junge \and Hun Hee Lee}

\address{Marius Junge :
Department of Mathematics
University of Illinois at Urbana-Champaign
273 Altgeld Hall 1409 W. Green Street Urbana, Illinois 61801, USA}
\email{junge@math.uiuc.edu}
\address{Hun Hee Lee : 
Department of Pure Mathematics
Faculty of Mathematics
University of Waterloo
200 University Avenue West, Waterloo, Ontario, Canada N2L 3G1}
\email{hh5lee@math.uwaterloo.ca, Lee.hunhee@gmail.com}


\begin{abstract}
We consider the $L_p$ norm estimates for homogeneous polynomials of $q$-gaussian variables ($-1\leq q\leq 1$).
When $-1<q<1$ the $L_p$ estimates for $1\leq p \leq 2$ are essentially the same as the free case ($q=0$),
whilst the $L_p$ estimates for $2\leq p \leq \infty$ show a strong $q$-dependence.
Moreover, the extremal cases $q = \pm 1$ produce decisively different formulae.

\end{abstract}

\maketitle

\section{Introduction}

In classical probability theory and in analysis the Wiener chaos is well understood and fundamental.
The orthgonal polynomials for the gaussian measure in $\Real^n$ are given by the Hermite polynomials.
In this paper we are interested in norm estimates for homogeneous polynomials in $q$-gaussian variables.
The family of $q$-gaussian variables, introduced by Bozejko and Speicher (\cite{BS91}),
is a natural noncommutative generalization of classical gaussian random variables.
Here $q$ ranges between $[-1,1]$. The case $q=1$ corresponds to the classical situation,
$q=-1$ reflects the fermionic case, and $q=0$ comes from Voiculescu's free probability theory.
Our work is motivated by the beautiful results by Haagerup/Pisier  (\cite{HP}) on the operator space structure of the space of homogeneous polynomials
in the free group. Similar result have been obtained by Ricard/Xu (\cite{RX}) for free products.
We want to show that norm estimates for polynomials of degree 2 and higher can detect the value $q$.

The space of polynomials of degree one is just the linear span of the ($q$-) gaussian variables and are, by now,
very well understood through Khintchine type inequalities.
The starting point of these result is Lust-Piquard's Khintchine inequality (\cite{LP}),
later extended by Lust-Piquard/Pisier to the range $1\le p\le 2$ (\cite{LPP}).
In the free case, Buchholz (\cite{Bu01}) provided very precise estimates.
For general information on free chaos we refer to the work of \cite{JPX-FreeChaos}, which stems from earlier work of Pisier/Parcet (\cite{ParPis05}).
For $p=\infty$, very nice estimates for polynomials can be found in \cite{RX}, the $L_p$ versions can be found in \cite{JPX-FreeChaos}.

Analyzing the known norm estimates for the span of the generators of $q$-gaussian random variables,
it turns out that the dependence on $q$ disappears.
Indeed, any suitable notion of gaussian random variables leads to similar expression for a fixed value of $p$.
In this paper we provide formulas for homogeneous polynomials of higher degree
and show that already for polynomials of degree $2$ the results do depend on $q$, however the dependence is quite subtle.

	\begin{thm}\label{11}
	Let $-1< q< 1$ and $(g_{q, i})_{i=1}^m$ (resp. $(g_i)^m_{i=1}$) be $q$-gaussian (resp. free-gaussian) random variables defined
	with respect to some reference state $\phi_q$ (resp. $\phi$) with density $D_q$ (resp. $D$).
	Let $(x_{ij})^m_{i,j=1}$ be $L_p(\M)$ valued coefficients for a von Neumann algebra $\M$ and $1\leq p\le 2$. Then
		$$ \norm{\sum_{i\neq j} x_{ij} \otimes D_q^{\frac{1}{2p}}g_{q, i} g_{q, j} D_q^{\frac{1}{2p}}}_p
		\sim_{c_q} \norm{\sum_{i\neq j} x_{ij} \otimes D^{\frac{1}{2p}}g_ig_jD^{\frac{1}{2p}}}_p.$$
	\end{thm}

The good news is that the right hand side can be calculated (up to some universal constant)
by a formula involving the best decomposition with respect to three norms (a generalization of the basic $\K$- functional spaces in interpolation).
Such decomposition norms should be considered classical in the theory as in so far
they appear already in Lust-Piquard/Pisier's operator valued Khintchine inequality for the range $1\le p\le 2$.
We refer to section \ref{sec-q-case} for a precise formulation. It seems that Theorem \ref{11} has no dependence in $q$ (except the constant).
However, there are singularities for $q=1$ and $q=-1$. In the specials cases we have to use a decoupling technique.
Indeed, it is well-known from Banach space theory (see e.g. the book of Ledoux and Talagrand (\cite{LT92})) that
	$$ \E \norm{\sum_{i\neq j} x_{ij} \otimes \eps_i\eps_j}_{X}\sim_c  \E \norm{\sum_{i\neq j} x_{ij} \otimes \eps_i \eps'_j}_X$$
holds for Banach space valued coefficients $x_{ij} \in X$.
Here $\eps_i$, $\eps_j'$ are independent coordinate functions defined on $\{-1, 1\}^n \times \{-1, 1\}^n$ and $\E$ is the corresponding expectation.
Based on new decoupling techniques we prove the following.

	\begin{thm}\label{qq} 
	Let $q=\pm 1$ and $(g_{q, i})_{i=1}^m$ be $q$-gaussian random variables defined
	with respect to some reference state $\phi_q$ with density $D_q$, and let $(g_i)^m_{i=1}$ and $D$ be as in Theorem \ref{11}.
	Let $1\le p<\infty $ and $(x_{ij})^m_{i,j=1}$ be $L_p(\M)$ valued coefficients for a von Neumann algebra $\M$.
 		\begin{enumerate}
	 		\item[i)] Let $q=-1$ and $x_{ij}= -x_{ji}$. Then
		  $$\norm{\sum_{i,j} x_{ij} \otimes D_{-1}^{\frac{1}{2p}}g_{-1,i} g_{-1,j} D_{-1}^{\frac{1}{2p}}}_p \sim
		  \norm{\sum_{i,j} x_{ij} \otimes D^{\frac{1}{2p}} g_i D^{\frac{1}{2p}} \otimes D^{\frac{1}{2p}} g_j D^{\frac{1}{2p}}}_p .$$
  
 \item[ii)]
   Let $q=1$ and $x_{ij}=x_{ji}$. Then
  $$\norm{\sum_{i,j} x_{ij} \otimes D_{1}^{\frac{1}{2p}}g_{1,i} g_{1,j} D_{1}^{\frac{1}{2p}}}_p \sim
		  \norm{\sum_{i,j} x_{ij} \otimes D^{\frac{1}{2p}}g_i D^{\frac{1}{2p}} \otimes D^{\frac{1}{2p}} g_j D^{\frac{1}{2p}}}_p .$$
   
\end{enumerate}
\end{thm}

Here again the right hand side can be calculated using martingale inequalities for the linear terms.
This leads to a four term maximum for $p\geq 2$ and a four term decomposition for $1\le p\le 2$.
The formulae are decisively different from the three term expression in Theorem \ref{11}.
Let us refer to section \ref{sec-CAR} for the concrete expressions.
Moreover, in section \ref{sec-q-case} we also extend this result to polynomials of arbitrary degree.
The notion becomes rather involved and the estimates depend on the degree of the polynomial.
In spirit our method is closely related to similar estimates for polynomials on the free group by Haagerup and Piser.
However, in our approach the decoupling is derived from Speicher's central limit procedure combined with the ultraproduct technique from \cite{J-Araki}.

Finally, let us come back to the $q$-dependence for polynomials of degree $2$.
For $p\geq 2$ we can use duality arguments starting from Theorem \ref{11}.
These estimates are based on the previous work of Nou (\cite{N04}).
Indeed, a classical tool for studying $q$-gaussian variables is the Wick order.
This order essentially implements the identification between the GNS-Hilbert space given by the vacuum state and the Fock space realization.
In our situation the usual product $g_{q, i} g_{q, j}$ ($i\neq j$) coincides with the Wick product,
however understanding $g_{q, i} g_{q, j}$ as a Wick product plays an essential role when $p\geq 2$.

	\begin{thm}
	Let $-1<q<1$, $(x_{ij})^m_{i,j =1}$ be $L_p(\M)$ valued coefficients for a von Neumann algebra $\M$ and $2\leq p \leq \infty$.
	Let $(g_{q,i})^m_{i=1}$, $D_q$, $(g_i)^m_{i=1}$ and $D$ be as in Theorem \ref{11}. Then
		$$\norm{\sum_{i \neq j} x_{ij}\otimes D_q^{\frac{1}{2p}}g_{q, i} g_{q, j} D_q^{\frac{1}{2p}}}_p
		\sim_{c_q} \norm{\sum_{i \neq j} (x_{ij} + q x_{ji}) \otimes D^{\frac{1}{2p}}g_i g_j D^{\frac{1}{2p}}}_p.$$
Moreover, the span of polynomials of degree $2$ is completely complemented in the corresponding $L_p$ space with a constant depending on $q$.
\end{thm}

The $q$-dependent term $x_{ij} + q x_{ji}$ above comes from the symmetrization operator $P_2$
on $\Hi^2$ for a Hilbert space $\Hi$ defined by
	$$P_2(f_1 \otimes f_2) = f_1\otimes f_2 + q f_2\otimes f_1$$
for any $f_1, f_2 \in \Hi$. See section \ref{subsec-Fock-Realization} for the details.
Again the result also holds for polynomials of higher degree, but we refer to the text for the explicit formulation.

The paper is organized as follows. In section \ref{sec-prelim} we present some preliminaries we need in sequel.
That includes the Fock space realization of the (generalized) $q$-gaussian random variables for $-1\leq q\leq 1$,
Wick product and some modular theory for the case $-1<q<1$, and matrix models for the case $q= \pm 1$.
At the end of the section \ref{sec-prelim} we collect some notations manly concerned with complicated indices we will encounter.
In section \ref{sec-q-case} we focus on the case $-1<q<1$. We first establish the free case using the result in \cite{JPX-FreeChaos}
and obtain an appropriate interpolation scale.
Using this interpolation scale and the approach of Nou (\cite{N04}) we can get the result for the general case $-1<q<1$.
In the final section we consider the remaining cases $q = \pm 1$.

We assume that the reader is familiar with standard concepts in operator algebra (\cite{Ta79, Ta03}), operator space theory (\cite{ER00, P03}),
noncommutative $L_p$ spaces (\cite{JX03, PX03}), and the related complex interpolation theory (\cite{Ko84, Te82, PX03}).
For a von Neumann algebra $\M$ we denote the noncommutative $L_p$ ($1\leq p \leq \infty$) space with respect to $\M$ by $L_p(\M)$.
When $E \subseteq L_p(\N)$ for another von Neumann algebra $\N$,
the norm closure of the algebraic tensor product $L_p(\M) \otimes E$ in $L_p(\M \bar{\otimes} \N)$ will be denoted by $L_p(\M;E)$.
For a Hilbert space $H$ we write $S_p(H) = L_p(B(H))$ and $H_r$ and $H_c$ imply the row and the column Hilbert space on $H$, respectively.
When $H = \ell_2$ we simply write $S_p$, $R$ and $C$, respectively.
$R_p$ and $C_p$ imply the linear space of the first row and column of $S_p$, respectively. 

\section{Preliminaries and Notations}\label{sec-prelim}

\subsection{The Fock space realization}\label{subsec-Fock-Realization}
We start with the Fock space realization of $q$-commutation relations for $-1\leq q\leq 1$.
Let $\Hi$ be an infinite dimensional separable complex Hilbert space
equipped with an orthonormal basis $(e_{\pm k})_{k\geq 1}$. We denote by $\F_0(\Hi)$ the associated free Fock space
$$\F_0(\Hi) = \mathbb{C}\Om \oplus \bigoplus_{n\geq 1}\Hi^{\otimes n},$$ where $\Om$ is a unit vector called vacuum.
We consider the operator of symmetrization $P_n$ on $\Hi^{\otimes n}$ defined by $$P_0 \Om = \Om,$$
$$P_n (f_1\otimes \cdots \otimes f_n) = \sum_{\pi \in S_n}q^{i(\pi)}f_{\pi(1)}\otimes \cdots \otimes f_{\pi(n)},$$
where $S_n$ denotes the symmetric group of permutations of $n$ elements and 
$$i(\pi) = \# \{(i,j) | 1\leq i,j \leq n, \pi(i) > \pi(j)\}$$ is the number of inversions of $\pi \in S_n$.

Now we define the {\it $q$-inner product} $\left\langle \cdot, \cdot \right\rangle_q$ on $\F_0(\Hi)$ by 
$$\left\langle \xi, \eta \right\rangle_q =\delta_{n,m}\left\langle \xi, P_n \eta \right\rangle
\;\, \text{for}\;\, \xi\in \Hi^{\otimes n}, \eta\in\Hi^{\otimes m},$$
where $\left\langle \cdot, \cdot \right\rangle$ is the inner product in $\Hi$.
Since $P_n$'s are strictly positive for $-1<q<1$ (\cite{BS91}), $\left\langle \cdot, \cdot \right\rangle_q$ is actually an inner product.
In this case, the Hilbert space $\F_0(\Hi)$ equipped with $q$-inner product is called {\it $q$-Fock space} on $\Hi$, and we denote it by $\F_q(\Hi)$.
When $q = \pm 1$, $P_n$'s are just positive, so we define $\F_q(\Hi)$ as the quotient of $(\F_0(\Hi), \left\langle \cdot, \cdot \right\rangle_q)$
by the corresponding kernel.

For two sequences of strictly positive reals $\lambda = (\lambda_k)_{k\geq 1}$ and $\mu = (\mu_k)_{k\geq 1}$
we define {\it $q$-(generalized) gaussian variables (or $q$-generalized circular variables)} by 
$$g_{q,k} = \lambda_k\ell_q(e_k) + \mu_k \ell^*_q(e_{-k}),$$ 
where $\ell_q(h)$ is the left creation operator by $h\in \Hi$ and $\ell^*_q(h)$ is the adjoint of $\ell_q(h)$.
It is easy to check that $g_{q,k}$'s satisfy the $q$-commutation relations
$$g^*_{q,k} \cdot g_{q, j} - q \cdot g_{q, j} \cdot g^*_{q,k} = \delta_{kj}(\lambda^2_k + \mu^2_k)I.$$
When $q= \pm 1$ we have additional relations $$g_{q,k} \cdot g_{q, j} - q \cdot g_{q, j} \cdot g_{q,k} = 0$$
which implies that $g_{q,k}$'s are CAR and CCR sequences in the corresponding cases.

We focus on $\Gamma_q$ ($-1\leq q <1$), the von Neumann algebra generated by $\{g_{q,k}\}_{k \geq 1}$.
When $q=1$ we define $\Gamma_1$ by the von Neumann algebra generated by
$\{\text{exp}(i \cdot g_{1,k})\}_{k \geq 1}$ since $g_{1,k}$'s are unbounded operators in this case.

There is a canonical way to translate the above picture into the framework of Shlyakhtenko and Hiai (\cite{S97, Hi}).
According to section 4 in \cite{S97} we can associate $g_{q,k}$ with an action
	$$U^k_t = \left( \begin{array}{cc}\cos\theta_k t & -\sin\theta_k t\\ \sin\theta_k t& \cos\theta_k t \end{array} \right),$$
where $\theta_k = \log [(\frac{\mu_k}{\lambda_k})^2]$ on $H_k (\cong \mathbb{R}^2)$.
Note that the basis $(\widetilde{e}_k, \widetilde{e}_{-k})$ on $H_k$ for this matrix representation is given by
\begin{equation}\label{basis-of-action}
\left(\begin{array}{cc}\widetilde{e}_k\\ \widetilde{e}_{-k} \end{array}\right)
= \frac{1}{\sqrt{\alpha_k +1}}\left( \begin{array}{cc} -i & \sqrt{\alpha_k}i \\ 1 & \sqrt{\alpha_k} \end{array} \right)
\left(\begin{array}{cc}e_k\\ e_{-k} \end{array}\right) = V\left(\begin{array}{cc}e_k\\ e_{-k} \end{array}\right),
\end{equation}
where $\alpha_k = \lambda^{-2}_k \mu^2_k$.

Let $s(\widetilde{e}_{k})$ and $s(\widetilde{e}_{-k})$ be semi-circular variables defined by
$$s(\widetilde{e}_{\pm k}) = \frac{1}{2}(\ell(\widetilde{e}_{\pm k}) + \ell^*(\widetilde{e}_{\pm k})).$$
Then we have
\begin{equation}\label{translation}
g_{q,k} = \frac{\sqrt{\lambda^2_k + \mu^2_k}}{2}(s(\widetilde{e}_k) + is(\widetilde{e}_{-k})).
\end{equation}
Since $\Hi$ is the complexification of $H_{\mathbb{R}} = \bigoplus_{k \geq 1}H_k$, by setting
$$U_t = \bigoplus_{k \geq 1}U^k_t$$ we get that $$\Gamma_q(H_{\mathbb{R}}, U_t) = \{s(e_{\pm k}) : k\geq 1\}'' = \Gamma_q.$$
In \cite{S97} we introduced another inner product $\left\langle \cdot, \cdot \right\rangle_U$ on $\Hi$ defined by
$$\left\langle x, y \right\rangle_U = \left\langle 2A(1+A)^{-1}x, y \right\rangle,$$ where $A$ is the operator satisfying $$U_t = A^{it}.$$
Note that $A = \bigoplus_{k \geq 1} A_k$ with $U^k_t = A_k^{it}$,
and by taking conjugate with respect to the basis change matrix $V$ in \eqref{basis-of-action}
we get the matrix representation of $A_k$ with respect to $(e_k, e_{-k})$ by
$$V^{-1}A_k V = \left( \begin{array}{cc}\lambda^2_k \mu^{-2}_k & 0 \\ 0 & \lambda^{-2}_k \mu^{2}_k \end{array} \right),$$
which implies that
\begin{align}\label{action-A}
A e_k = \lambda^{2}_k \mu^{-2}_k e_k
\end{align}
for any $k\geq 1$.

\subsection{The Wick product and some modular theory for the case $-1<q<1$}

Let $-1<q<1$. Since it is well known that $\Om$ is separating for $\Gamma_q$, for every $\xi \in \Gamma_q \Om$ there exist a unique operator
$W(\xi) \in \Gamma_q$ such that $$W(\xi)\Om = \xi.$$ This $W$ is called the {\it Wick product}.
There is a useful decomposition of the Wick product as follows.
\begin{equation}\label{eq-Wick}
W(\xi) = \sum^n_{k=0}U_k R^*_{n,k}(\xi)
\end{equation}
for $\xi \in \Hi^{\otimes n}$, where $R_{n,k}$ is the operator on $\Hi^{\otimes n}$ given by
$$R_{n,k}(f_1\otimes \cdots \otimes f_n) = \sum_{\pi \in S_n / S_{n-k} \times S_{k}}q^{i(\pi)}f_{\pi^{-1}(1)}\otimes \cdots \otimes f_{\pi^{-1}(n)}$$
and $U_k : \Hi^{\otimes n-k}_c \otimes \Hi^{\otimes k}_r \rightarrow B(\F_q(\Hi))$ is the operator defined by
$$U_k (e_1\otimes \cdots \otimes e_n) = \ell(e_1)\cdots \ell(e_{n-k})\ell^*(e_{n-k+1})\cdots \ell^*(e_n).$$
In the above $S_n / S_{n-k} \times S_{k}$ means the representatives of the right cosets of $S_{n-k} \times S_{k}$ in $S_n$
with minimal numbers of inversions.
It is well known that
	\begin{equation}\label{P-n-R-nk}
	P_n = R_{n,k}(P_{n-k} \otimes P_{k})
	\end{equation}
and by Corollary 1 in \cite{N04}
	\begin{equation}\label{U-k}
	\norm{U_k : H^{\otimes n-k}_c \otimes_h H^{\otimes k}_r \rightarrow B(\F_q(\Hi))}_{cb} \leq C_q,
	\end{equation}
where $H^{\otimes k}$ is the Hilbert space on $\Hi^{\otimes k}$ equipped with the $q$-inner product and
	$$C_q = \prod_{n \geq 1}\frac{1}{1-q^n}.$$

We can estimate the norm of $P_n$ and its inverse. In \cite{BS94} it is shown that
	$$P_n \leq \frac{1}{1-q} I_1 \otimes P_{n-1},$$ where $I_1$ is the formal identity on $H$, thus by repeating the above operator inequality we have
	\begin{equation}\label{norm-P-n}
	\norm{P_n : \Hi^{\otimes n} \rightarrow \Hi^{\otimes n}} \leq \Big(\frac{1}{1-q}\Big)^{n-1}.
	\end{equation}
Moreover, it is shown in \cite{B99} that $I_1\otimes P_{n-1} \leq w(q)^{-1}P_n$ for some constant $w(q)>0$, thus by repeating again we have
	\begin{equation}\label{norm-P-n-inverse}
	\norm{P^{-1}_n : \Hi^{\otimes n} \rightarrow \Hi^{\otimes n}} \leq w(q)^{-n+1}.
	\end{equation}

From the definition it is clear that
	\begin{equation}\label{P-n-isometry}
	P^{\frac{1}{2}}_n : H^{\otimes n} \rightarrow \Hi^{\otimes n}
	\end{equation}
is an isometry.

We close this section with some modular theory for $\Gamma_q$ and $\phi_q$, the vacuum state defined by
$\phi_q(\cdot) = \left\langle \Om \;\cdot, \Om \right\rangle_q$.
It is well known that the modular group $\sigma_t$ with respect to $\phi_q$ satisfies the following.
$$\sigma_t(g_{q,k}) = (\lambda^{-1}_k \mu_k)^{2it}g_{q,k}.$$
Thus, $g_{q,k}$ is an analytic element satisfying
\begin{equation}\label{gq-in-Lp}
D^{\frac{1}{2p}}_q g_{q,k}  = (\lambda^{-1}_k \mu_k)^{\frac{1}{p}} g_{q,k} D^{\frac{1}{2p}}_q.
\end{equation}
Recall that the anti-linear map $S$ is the closure of the operator given by $$S(x\Om) = x^*\Om$$ for all $x\in \Gamma_q$.
Then $S$ can be written as $S = J\Delta^{\frac{1}{2}}$, where $J$ is the modular conjugation and $\Delta$ is the modular operator.
By \cite{Hi} we know that $\Delta$ is the closure of the operator $\bigoplus_{n\geq 0}(A^{-1})^{\otimes n}$, which implies
\begin{equation}\label{Delta-info}
S|_{H^{\otimes n}} = J(A^{-\frac{1}{2}})^{\otimes n}.
\end{equation}

\subsection{Matrix models for CAR and CCR generators}\label{subsec-matrix}

Now we focus on the case $q = \pm 1$. We will follow the approach in Example 3.8 of \cite{J-Araki} to construct matrix model for $g_{q,k}$'s.
First, we fix $m \in \n$ and restrict our attention to $g_{q, 1}, \cdots g_{q, m}$.
By \eqref{translation} we need $2m$ semi-circular variables to approximate $g_{q,k}$'s.
Thus, we consider an associated weight $\psi$ on $N = \ell^{2m}_{\infty}(M_2)$.
Let $K$ be the real Hilbert space consisting of self-adjoint elements of $N$ equipped with the inner product
$$\left\langle x,y \right\rangle_K = \frac{\psi(xy) + \psi(yx)}{2}.$$
Then $\widetilde{U}_t(m) = U_t(m) \otimes I$ is a one-parameter group of unitaries on $K$, where $$U_t(m) = \bigoplus^m_{k=1}U^k_t.$$

By the relationship $$\psi(xy) = \left\langle x, y \right\rangle_{\widetilde{U}(m)}$$ for self-adjoint elements $x,y$
in Example 3.8 of \cite{J-Araki} we can determine $\psi$ as follows.
$$\psi(x) = \sum_{1\leq \abs{k}\leq m}[(2-\sigma_k)x_{11}(k) + \sigma_k x_{22}(k)]$$
for $x = [(x_{ij}(\pm 1)), \cdots, (x_{ij}(\pm m))] \in \ell^{2m}_{\infty}(M_2)$, where
$$\sigma_k = \frac{2\mu^2_{\abs{k}}}{\lambda^2_{\abs{k}} + \mu^2_{\abs{k}}}.$$ Let
\begin{equation}\label{def-u_n}
u_n(x) = \sqrt{\frac{4m}{n}}\sum^n_{k=1}v_k \otimes \pi_k(x),
\end{equation}
where $\pi_k : N \rightarrow N^{\otimes n}$ is the homomorphism which sends $N$ in the $k$-th component
and $v_k\in M_{2^n}$ are self-adjoint unitaries such that
\begin{equation}\label{relation-v_k}
\text{$v_k v_j = - v_jv_k$ (when $q=-1$) or $v_k v_j = v_jv_k$ (when $q=1$)}.
\end{equation}
Note that we are using the scaling factor $4m$ since $\psi(1_N) = 4m$ in this situation.
Actually, we will use the following special choice of $v_k$'s.
Let $$v_{k,k} = \left( \begin{array}{cc} 0 & 1\\ 1 & 0\end{array} \right)\;\, \text{and for}\;\, j<k \;\,
v_{j,k} = \left( \begin{array}{cc} 1 & 0\\ 0 & q\end{array} \right).$$ Then we set
\begin{equation}\label{def-v_k}
v_k = v_{1,k}\otimes \cdots \otimes v_{k,k}\otimes 1\otimes \cdots \otimes 1\in M_{2^n}.
\end{equation}

Let $\U$ be a free ultrafilter on $\n$ and $$N_\U = \Big(\prod_{n,\U} (M_{2^n} \otimes N^{\otimes n})_*\Big)^*$$ with the ultraproduct state
$$\phi_\U = (\tau_n\otimes \Big(\frac{\psi}{4m}\Big)^{\otimes n})_{n,\U},$$ where $\tau_n$ is the normalized trace on $M_{2^n}$.
Let $(\delta_k)_{1\leq \abs{k} \leq m}$ be the unit vectors  in $\ell^{2m}_\infty$.
In \cite{J-Araki} it is shown that $u_n(\delta_k \otimes e_{12})$ (with respect to $\phi_\U$) converges in $*$-distribution to
$s(\widetilde{e}_k)$ (with respect to $\phi$) as $n$ goes infinity,
and if we consider the map $$\Phi : L_p(\Gamma_{\pm 1})  \rightarrow L_p(N_\U)\;\, (1\leq p < \infty)$$ defined by
$$\Phi(D^{\frac{1}{2p}}_{\pm 1} P(s(\widetilde{e})_{\pm 1}, \cdots, s(\widetilde{e})_{\pm m})D^{\frac{1}{2p}}_{\pm 1})
= D^{\frac{1}{2p}}_{\phi_\U} P(u_n(\delta_{\pm 1} \otimes e_{12}), \cdots, u_n(\delta_{\pm m} \otimes e_{12}))D^{\frac{1}{2p}}_{\phi_\U}$$
for any noncommutative polynomial $P$, then we know that
\begin{equation}\label{matrix-model-inclusion}
I_\M \otimes \Phi : L_p(\M \bar{\otimes}\Gamma_{\pm 1}) \rightarrow L_p(\M \bar{\otimes}N_\U)
\end{equation}
is an isometry for any von Neumann algebra $\M$.

When $q=1$ we shall understand $s(\widetilde{e}_k) D^{\frac{1}{2p}}_{1}$ as
$\frac{d}{i dt}(\text{exp}(it s(\widetilde{e}_k))D^{\frac{1}{2p}}_{1})|_{t=0}$.

\subsection{Notations}

In the following we will frequently use the index sets $$\I^d = \{\und{i} = (i_1, \cdots, i_d) \in \n^d : i_1 \neq \cdots \neq i_d\}$$
and for $m\in \n$ $$\I^d_m = \{\und{i} = (i_1, \cdots, i_d) \in \n^d : 1\leq i_1 \neq \cdots \neq i_d \leq m\},$$
where $i_1 \neq \cdots \neq i_d$ means that there is no repeating index. 
If we remove the restriction of repetition, then we have the index sets $\n^d$ and
$$\n^d_m = \{\und{i} = (i_1, \cdots, i_d) \in \n^d : 1\leq i_1,\cdots,i_d \leq m\},$$ respectively 

For $\und{i}\in \I^d_m \; \text{or}\; \n^d_m$ we will use the notation
	$$g_{q, \und{i}} := g_{q,i_1}\cdots g_{q,i_d},$$
	$$g_{p, q,\und{i}} := D_q^{\frac{1}{2p}}g_{q,i_1}\cdots g_{q,i_d}D_q^{\frac{1}{2p}}$$
and
	$$e_{\und{i}} = e_{i_1}\otimes \cdots \otimes e_{i_d}.$$
Similarly, for a sequence of real numbers $\lambda = (\lambda_k)_{k \geq 1}$ we will write
	$$\lambda_{\und{i}} = \lambda_{i_1}\cdots \lambda_{i_d}.$$

We use the symbol $a\lesssim b$ if there is a $C>0$ such that $a \leq C b$ and $a\sim b$ if $a\lesssim b$ and $b\lesssim a$.
Similarly, $a\lesssim_{c_d} b$ (resp. $a\lesssim_{c_{q, d}} b$) if there is a constants $C>0$ depending only on $d$ (resp. only on $q$ and $d$)
such that $a \leq C b$. We write $a\sim_{c_d} b$ if $a\lesssim_{c_d} b$ and $b\lesssim_{c_d} a$, and the meaning of $a\sim_{c_{q,d}} b$ is similar.

From now on we fix two sequences of strictly positive reals $\lambda = (\lambda_k)_{k \geq 1}$ and $\mu = (\mu_k)_{k \geq 1}$,
the number of degree $d \in \n$ and a von Neumann algebra $\M$.
$\M$ is equipped with a distinguished normal semifinite faithful weight $\varphi$ which induces the trace functional $\text{tr}_\M$ on $L_1(\M)$.
We denote the unit of $\M$ by $1_\M$ whilst $I_\M$ implies the identity map on $\M$.

\section{The case $-1 < q < 1$}\label{sec-q-case}

\subsection{Free case: Building a model}
First, we consider the free case ($q=0$) and obtain an interpolation result as a corollary.
We will simply write $g_{k}$, $\Gamma$ and $\phi$ instead of $g_{0,k}$, $\Gamma_0$ and $\phi_0$, respectively.
For $\und{i}\in \I^d_m \; \text{or}\; \n^d_m$ we will use the notations
	$$g_{\und{i}} := g_{i_1}\cdots g_{i_d}\;\,\text{and}\;\,g_{p,\und{i}} := D^{\frac{1}{2p}}g_{i_1}\cdots g_{i_d}D^{\frac{1}{2p}}.$$

For $x = (x_{\und{i}})_{\und{i}\in \I^d_m} \subseteq L_p(\M)$ and $0\leq k\leq d$ we denote by
	$$\R\C_p^{d,k}(\lambda, \mu; x) = \norm{\sum_{\und{i}\in \I^d_m}x_{\und{i}}\otimes
	f^{p, i_1}\otimes \cdots \otimes f^{p, i_{k}} \otimes f_{p, i_{k+1}} \otimes \cdots \otimes f_{p, i_d}}_{L_p(\M ; E^{d,k}_p)},$$
where
	$$f^{p,i} = \lambda^{\frac{1}{p'}}_i \mu^{\frac{1}{p}}_i e_{i1} \in C_p,\;\, f_{p,i} = \lambda^{\frac{1}{p}}_i \mu^{\frac{1}{p'}}_i e_{1i}\in R_p$$
and
	$$E^{d,k}_p = C^{\otimes k}_p \otimes_h R^{\otimes d-k}_p \subseteq S_p(\ell^{\otimes d}_2).$$
We define the corresponding $\K$- and $J$-functional spaces $\K^d_p(\lambda, \mu)$ $(1\leq p\leq 2)$ and $\J^d_p(\lambda, \mu)$ $(2\leq p\leq \infty)$
as the closures of finite tuples in $L_p(\M)$ indexed by $\I^d$ with respect to the following norms. 
$$\norm{x}_{\K^d_p(\lambda,\mu)}=\inf\Big\{\R\C_p^{d,0}(\lambda,\mu;x^0)+\R\C_p^{d,1}(\lambda,\mu;x^1)+\cdots +\R\C_p^{d,d}(\lambda,\mu; x^d)\Big\},$$
where the infimum runs over all possible $x_{\und{i}} = x^0_{\und{i}} + \cdots + x^d_{\und{i}}$ and
$x^k = (x^k_{\und{i}})_{\und{i}\in \I^d_m}$ for $0\leq k \leq d$, and
$$\norm{x}_{\J^d_p(\lambda, \mu)}
= \max_{0\leq k \leq d}\Big\{ \R\C^{d,0}_p(\lambda, \mu; x),\; \R\C^{d,1}(\lambda,\mu; x)_p,\;\cdots,\;\R\C^{d,d}_p(\lambda, \mu; x)\Big\}.$$
Note that $\J^d_{p'}(\lambda,\mu) = (\K^d_{p}(\lambda,\mu))^*$ with the duality bracket
$$\left\langle (x_{\und{i}})_{\und{i}\in \I^d}, (z_{\und{i}})_{\und{i}\in \I^d} \right\rangle
= \sum_{\und{i}\in \I^d} \lambda_{\und{i}}\mu_{\und{i}}\text{tr}_{\M}(x^*_{\und{i}}z_{\und{i}}).$$

\begin{thm}\label{thm-free}
Let $m\in \n$ and $x = (x_{\und{i}})_{\und{i}\in \I^d_m} \subseteq L_p(\M).$
Then we have the following equivalence.
	\begin{align*}
	\norm{\sum_{\und{i}\in \I^d_m}x_{\und{i}}\otimes g_{p,\und{i}}}_{L_p(\M \bar{\otimes}\Gamma)} \sim_{c_d} 		
	\left\{\begin{array}{ll}\norm{x}_{\K^d_p(\lambda, \mu)} & \text{for\, $1\leq p \leq 2$}\\
	\norm{x}_{\J^d_p(\lambda, \mu)} & \text{for\, $2\leq p \leq \infty$.}\end{array} \right.
	\end{align*}
Moreover, $\G^d_p$, the closed subspace of $L_p(\Gamma)$ spanned by $$\{g_{p,\und{i}} : \und{i}\in \I^d\}$$
is completely complemented with the constant depending only on $d$. 
\end{thm}
\begin{proof}
We consider the following maps.
\begin{align*}
u_p : & \left\{\begin{array}{l}\K^d_p(\lambda, \mu)\; (1\leq p\leq 2)\\ \J^d_p(\lambda, \mu)\; (2\leq p \leq \infty) \end{array} \right.
 \rightarrow \;\;\;\; L_p(\M \bar{\otimes} \Gamma),\\
& \;\;\;\;\;\;\;\;(x_{\und{i}})_{\und{i}\in \I^d_m} \;\;\;\;\;\;\;
\mapsto \sum_{\und{i}\in \I^d_m}x_{\und{i}}\otimes g_{p,\und{i}}.
\end{align*}
What we need to do is to show that $u_p$'s are isomorphisms with bounded constants.
By a usual density argument it is enough to check that $u_p$'s are isomorphisms
when they are restricted to the span of $(x_{\und{i}})_{\und{i}\in \I^d_m}$ for a fixed $m\in \N$ with constants independent of $m$.

First, we consider the case $2\leq p \leq \infty$ and will use an induction on $d$.
When $d=1$ we are done by Theorem 3.1 of \cite{X-Lp-Grothendieck} (or Theorem 5.1 of \cite{JPX-FreeChaos}).
Suppose we have the result for $d$ and consider the case for $d+1$.

We will simply write $(i, \und{i}) \in \I^{d+1}_m$ when we have
$$\text{$\und{i}\in \I^d_m$ and $1\leq i \leq m$ satisfying $i \neq i_k$, $1\leq k\leq d$.}$$
In this case we will use the notation $x_{i, \und{i}}$ and $g_{i, \und{i}}$ instead of $x_{\und{i}'}$ and $g_{\und{i}'}$, respectively,
where $$\und{i}' = (i, i_1, \cdots, i_d) \in \I^{d+1}_m.$$

For $x = (x_{i, \und{i}})_{(i, \und{i}) \in \I^{d+1}_m} \subseteq L_p(\M)$ we consider
\begin{align*}
\sum_{(i, \und{i}) \in \I^{d+1}_m}x_{i, \und{i}}\otimes D^{\frac{1}{2p}}g_{i, \und{i}} D^{\frac{1}{2p}}
& = \sum_{(i, \und{i}) \in \I^{d+1}_m}(1_{\M}\otimes D^{\frac{1}{2p}}g_i)(x_{i, \und{i}} \otimes g_{\und{i}} D^{\frac{1}{2p}})\\
& = \sum_{(i, \und{i}) \in \I^{d+1}_m}\alpha_{i}\beta_{i, \und{i}}.
\end{align*}
By applying Theorem B of \cite{JPX-FreeChaos} we have
\begin{align*}
\lefteqn{\norm{\sum_{(i, \und{i}) \in \I^{d+1}_m}\alpha_{i}\beta_{i,\und{i}}}_p}\\
& \sim \norm{\bigg(\sum_{(i, \und{i}), (j, \und{j}) \in \I^{d+1}_m}\beta^*_{i, \und{i}}
\E(\alpha^*_{\und{i}}\alpha_{j})\beta_{j, \und{j}}\bigg)^{\frac{1}{2}}}_p
+ \norm{\bigg(\sum_{(i, \und{i}), (j, \und{j}) \in \I^{d+1}_m}\alpha_{i}
\E(\beta_{i, \und{i}}\beta^*_{j, \und{j}})\alpha^*_{\und{j}}\bigg)^{\frac{1}{2}}}_p\\
& = A + B,
\end{align*}
where $\E : L_p(\M \bar{\otimes} \Gamma) \rightarrow L_p(\M)$ is the projection induced from the conditional expectation $I_{\M} \otimes \phi$ (see section 2 of \cite{JX03}). 

Since we have by \eqref{gq-in-Lp}
\begin{align*}
\E(\alpha^*_{\und{i}}\alpha_{j})
& = 1_\M \otimes \phi(g^*_i D^{\frac{1}{p}} g_j)\\
& = \lambda^{-\frac{2}{p}}_i\mu^{\frac{2}{p}}_j 1_\M \otimes \phi(g^*_i g_j)D^{\frac{1}{p}}\\
& = \delta_{ij}\lambda^{\frac{2}{p'}}_i\mu^{\frac{2}{p}}_j 1_\M \otimes D^{\frac{1}{p}}
\end{align*}
we get the following by the induction hypothesis.
\begin{align*}
A & = \norm{\bigg( \sum^n_{i=1} \lambda^{\frac{2}{p'}}_i\mu^{\frac{2}{p}}_i
(\sum_{(i, \und{i}), (i, \und{j}) \in \I^{d+1}_m} x^*_{i, \und{i}} x_{i, \und{j}} \otimes g^*_{p,\und{i}} g_{p,\und{j}})\bigg)^{\frac{1}{2}}}_p\\
& = \norm{\sum_{(i, \und{i}) \in \I^{d+1}_m} x_{i, \und{i}} \otimes g_{p,\und{i}} \otimes f^{p,i}}_p\\
& \sim_{c_d} \max\Bigg\{\norm{\sum_{(i, \und{i}) \in \I^{d+1}_m} x_{i, \und{i}}\otimes f^{p,i_1}\otimes \cdots \otimes f^{p,i_d}\otimes f^{p,i}}_p,\\
& \;\;\;\;\;\;\;\;\;\;\;\;\;\;\,
\norm{\sum_{(i, \und{i}) \in \I^{d+1}_m} x_{i, \und{i}}\otimes f_{p,i_1}\otimes f^{p,i_2}\otimes \cdots \otimes f^{p,i_d}\otimes f^{p,i}}_p, \;\cdots, \\
& \;\;\;\;\;\;\;\;\;\;\;\;\;\;\,
\norm{\sum_{(i, \und{i}) \in \I^{d+1}_m} x_{i, \und{i}}\otimes f_{p, i_1}\otimes \cdots \otimes f_{p, i_d}\otimes f^{p,i}}_p \Bigg\}.
\end{align*}
Similarly we have
\begin{align*}
\E(\beta_{i, \und{i}}\beta^*_{j, \und{j}}) & = (I_\M\otimes \phi)(x^*_{i, \und{i}} x_{j, \und{j}} \otimes g_{\und{i}}D^{\frac{1}{p}}g^*_{\und{j}})
= \lambda^{\frac{2}{p}}_{\und{i}}\mu^{-\frac{2}{p}}_{\und{j}}x^*_{i, \und{i}} x_{j, \und{j}} \otimes \phi(g_{\und{i}}g^*_{\und{j}})D^{\frac{1}{p}}\\
& = \delta_{i_1 j_1}\cdots \delta_{i_d j_d}\lambda^{\frac{2}{p}}_{\und{i}}\mu^{\frac{2}{p'}}_{\und{j}}x^*_{i, \und{i}} x_{j, \und{j}}
\otimes D^{\frac{1}{p}},
\end{align*}
so that
\begin{align*}
B & = \norm{\bigg( \sum_{\und{i}\in \I^d_m} \lambda^{\frac{2}{p}}_{\und{i}}\mu^{\frac{2}{p'}}_{\und{j}}
(\sum_{\substack{i,j \neq i_k\\ 0\leq k\leq d}} x_{i, \und{i}} x^*_{j, \und{j}}\otimes g_{p,i}g^*_{p,j})\bigg)^{\frac{1}{2}}}_p\\
& = \norm{\sum_{(i,\und{i}) \in \I^{d+1}_m} x_{i, \und{i}}
\otimes g_{p,i} \otimes f_{p,i_1}\otimes \cdots \otimes f_{p,i_d}}_p\\
& \sim \max\Bigg\{\norm{\sum_{(i,\und{i}) \in \I^{d+1}_m} x_{i, \und{i}}\otimes f_{p,i_1}\otimes \cdots \otimes f_{p,i_d}\otimes f_{p,i}}_p,\\
& \;\;\;\;\;\;\;\;\;\;\;\;\;\;\;
\norm{\sum_{(i,\und{i}) \in \I^{d+1}_m} x_{i, \und{i}}\otimes f_{p,i_1}\otimes \cdots \otimes f_{p,i_d}\otimes f^{p,i}}_p \Bigg\}.
\end{align*}
Consequently, we get the result for $d+1$.

Now we apply a duality argument for the case $1\leq p \leq 2$.
By the above result we know that $u_p$'s $(2\leq p \leq \infty)$ are isomorphisms with bounded constants independent of $m$.
Note that for $x = (x_{\und{i}})_{\und{i}\in \I^d_m}$, $z = (z_{\und{i}})_{\und{i}\in \I^d_m}$ and $1\leq p\leq 2$ we have by \eqref{gq-in-Lp}
	\begin{align}\label{duality1}
	\lefteqn{\left\langle u_p(x), u_{p'}(z) \right\rangle}\\
	& = \sum_{\und{i}, \und{j}\in \I^d_m}\text{tr}_{\M}(x^*_{\und{i}}z_{\und{j}}) \text{tr}_{\Gamma}(g^*_{p,\und{i}} g_{p',\und{j}}) \nonumber\\
	& =
	\sum_{\und{i},\und{j}\in\I^d_m}\lambda^{-\frac{1}{p}}_{\und{i}}\mu^{\frac{1}{p}}_{\und{i}}\lambda^{-\frac{1}{p'}}_{\und{j}}\mu^{\frac{1}{p'}}_{\und{j}}
	\text{tr}_{\M}(x^*_{\und{i}}z_{\und{j}}) \text{tr}_{\Gamma}(D^{\frac{1}{p}}g^*_{\und{i}} g_{\und{j}}D^{\frac{1}{p'}}) \nonumber\\
	& =
	\sum_{\und{i},\und{j}\in\I^d_m}\lambda^{-\frac{1}{p}}_{\und{i}}\mu^{\frac{1}{p}}_{\und{i}}\lambda^{-\frac{1}{p'}}_{\und{j}}\mu^{\frac{1}{p'}}_{\und{j}}
	\text{tr}_{\M}(x^*_{\und{i}}z_{\und{j}}) \phi(g^*_{\und{i}} g_{\und{j}}) \nonumber\\
	& = \sum_{\und{i}\in \I^d_m} \lambda_{\und{i}}\mu_{\und{i}}\text{tr}_{\M}(x^*_{\und{i}}z_{\und{i}})
	= \left\langle x, z \right\rangle. \nonumber
	\end{align}
Moreover, for any $Y \in L_p(\M \bar{\otimes} \Gamma)$ we have
	\begin{align}\label{duality2}
	\left\langle u_p(x), Y \right\rangle
	= \left\langle u_p(x), u_{p'} u^*_p(Y) \right\rangle.
	\end{align}

Combining \eqref{duality1} and \eqref{duality2} it is enough to show that
$u_{p'}u^*_p$'s $(1\leq p \leq 2)$ are bounded with constants independent of $m$.
It is straightforward to check that the maps
	$$u_{p'}u^*_p : L_{p'}(\M \bar{\otimes} \Gamma) \rightarrow L_{p'}(\M \bar{\otimes} \Gamma)$$
are essentially the same maps for $1\le p\le 2$ in the sense of interpolation theory.
Thus, we are only to check the norms of $u_{\infty}u^*_1$ and $u_2 u^*_2$ by complex interpolation.
Since $p=2$ case is trivial, it suffices to show that $u_1$ is bounded.

Note that $$L^r_2(\Gamma)\widehat{\otimes} L^c_2(\Gamma) \rightarrow L_1(\Gamma),\;\, a\otimes b \mapsto ab$$ is a complete contraction,
so that for $0\leq k\leq d$ we have
\begin{align*}
\lefteqn{\norm{\sum_{\und{i}\in \I^d_m}x_{\und{i}}\otimes D^{\frac{1}{2}}g_{\und{i}}D^{\frac{1}{2}}}_{L_1(\M \bar{\otimes}\Gamma)}}\\
& \leq \norm{\sum_{\und{i}\in \I^d_m}x_{\und{i}}\otimes D^{\frac{1}{2}}g_{i_1}\cdots g_{i_{k}}
\otimes g_{i_{k+1}}\cdots g_{i_d}D^{\frac{1}{2}}}_{L_1(\M)\widehat{\otimes} L^r_2(\Gamma) \widehat{\otimes} L^c_2(\Gamma)}.
\end{align*}
Let $\und{i}^k = (i_1, \cdots, i_k) \in \I^k_m$ and $\und{i}^{d-k} = (i_{k+1}, \cdots, i_d) \in \I^{d-k}_m$.
Since $\{D^{\frac{1}{2}}g_{\und{i}^k}\}_{\und{i}^k \in \I^k_m}$ and $\{g_{\und{i}^{d-k}}D^{\frac{1}{2}}\}_{\und{i}^{d-k} \in \I^{d-k}_m}$
are orthogonal families of vectors in $L_2(\Gamma)$ and
	$$\text{tr}_{\Gamma}(D^{\frac{1}{2}}g_{\und{i}^k}g^*_{\und{i}^k} D^{\frac{1}{2}})
	= \phi(g_{\und{i}^k}g^*_{\und{i}^k}) = (\mu_{i_1}\cdots \mu_{i_{k}})^2$$
and
	$$\text{tr}_{\Gamma}(D^{\frac{1}{2}}g^*_{\und{i}^{d-k}}g_{\und{i}^{d-k}}D^{\frac{1}{2}})
	= (\lambda_{i_{k+1}}\cdots \lambda_{i_{d}})^2,$$
we obtain complete isometries
	$$\text{span}\{D^{\frac{1}{2}}g_{\und{i}^k}\}_{\und{i}^k \in \I^k_n} (\subseteq L^r_2(\Gamma))
	\cong \text{span}\{f^{1,i_1}\otimes \cdots \otimes f^{1,i_{k}}\}_{\und{i}^k \in \I^k_n} \subseteq C^{\otimes k}_1$$
and
	$$\text{span}\{g_{\und{i}^{d-k}}D^{\frac{1}{2}}\}_{\und{i}^{d-k} \in \I^{d-k}_n} (\subseteq L^c_2(\Gamma))
	\cong \text{span}\{f_{1,i_{k+1}}\otimes \cdots \otimes f_{1,i_d}\}_{\und{i}^{d-k} \in \I^{d-k}_n} \subseteq R^{\otimes d-k}_1.$$
This implies $u_1$ is a contraction.

Now we have isomorphisms $u_p$ for $1\leq p \leq \infty$ with bounded constants (depending only on $d$).
Thus $\G^d_p$'s are completely complemented with constants depending only on $d$ by the following projections.
	$$u_{p}u^*_{p'} : L_p(\M \bar{\otimes} \Gamma) \rightarrow L_p(\M \bar{\otimes} \Gamma).$$

\end{proof}

\begin{rem}{\rm
The proof of Theorem \ref{thm-free} follows closely the Theorem F in \cite{JPX-FreeChaos},
but we have to go a step further in proving the complementation result.
}
\end{rem}

Now we get an interpolation scale we need in the following section.

\begin{cor}\label{cor-interpolation1}
$\{ \K^d_p(\lambda, \mu) : 1\leq p\leq 2\} \cup \{ \J^d_p(\lambda, \mu) : 2\leq p\leq \infty\}$
forms an interpolation scale. In particular, for $1< p <\infty$ and $\theta = \frac{1}{p}$ we have
	$$[\J^d_{\infty}(\lambda, \mu), \K^d_1(\lambda, \mu)]_{\theta} \cong
	\left\{ \begin{array}{ll} \K^d_p(\lambda, \mu) & \text{for $1\leq p\leq 2$}\\
	\J^d_p(\lambda, \mu) & \text{for $2\leq p\leq \infty$} \end{array}\right.$$
isomorphically with constants depending only on $d$.
\end{cor}

If we rephrase the above in the operator space language (see \cite{P98}), then we get the following interpolation result as a corollary.
Let $1 \leq p \leq \infty$ and $f^{p,i}$ and $f_{p,i}$'s are as before.
For $0\leq k \leq d$ we denote the closed subspace of $C^{\otimes k}_p \otimes_h R^{\otimes d-k}_p$ spanned by
$$\{f^{p, i_1}\otimes \cdots \otimes f^{p, i_k} \otimes f_{p, i_{k+1}} \otimes \cdots \otimes f_{p, i_d}: \und{i} \in \I^d\}$$
by $RC^{d,k}_p(\lambda, \mu)$. Now we define
	$$RC^d_p(\lambda, \mu) = \left\{ \begin{array}{ll} RC^{d,0}_p(\lambda, \mu) + RC^{d,1}_p(\lambda, \mu) + \cdots + RC^{d,d}_p(\lambda, \mu)
	& \text{for $1\leq p\leq 2$}\\
	RC^{d,0}_p(\lambda, \mu) \cap RC^{d,1}_p(\lambda, \mu) \cap \cdots \cap RC^{d,d}_p(\lambda, \mu)
	& \text{for $2\leq p\leq \infty$} \end{array}\right..$$

\begin{cor}\label{cor-interpolation2}
$\{RC^d_p(\lambda, \mu) : 1\leq p\leq \infty\}$ forms an interpolation scale, i.e. for $0< \theta <1$ and $1\leq p, p_0, p_1 \leq \infty$
satisfying $\frac{1-\theta}{p_0} + \frac{\theta}{p_2} = \frac{1}{p}$ we have
	$$[RC^d_{p_1}(\lambda, \mu), RC^d_{p_2}(\lambda, \mu)]_{\theta} \cong RC^d_p(\lambda, \mu)$$
completely isomorphically with constants depending only on $d$.
\end{cor}

\subsection{The case of general $q$}

For $x = (x_{\und{i}})_{\und{i}\in \I^d_m} \subseteq L_p(\M)$ we can write
	$$x = \sum_{\und{i}\in \I^d_m} x_{\und{i}} \otimes e_{\und{i}}$$
as an element of $\J^d_p(\lambda,\mu)$.
Since $P_d$ is defined on $\Hi^{\otimes d}$ the following is well-defined.
	$$(I_{L_p(\M)}\otimes P_d)x := \sum_{\und{i}\in \I^d_m} x_{\und{i}} \otimes P_d(e_{\und{i}}).$$
It is clear that $(I_{L_p(\M)}\otimes P_d)x$ is also an element of $\J^d_p(\lambda,\mu)$.

\begin{thm}\label{thm-q}
Let $m\in \n$ and $x = (x_{\und{i}})_{\und{i}\in \I^d_m} \subseteq L_p(\M).$ Then we have the following equivalence.
\begin{align*}
\norm{\sum_{\und{i}\in \I^d_m}x_{\und{i}}\otimes g_{p,q,\und{i}}}_{L_p(\M \bar{\otimes}\Gamma_q)} \sim_{c_{d,q}} \left\{\begin{array}{ll}\norm{x}_{\K^d_p(\lambda, \mu)} & \text{for\, $1\leq p \leq 2$}\\
\norm{(I_{L_p(\M)}\otimes P_d)x}_{\J^d_p(\lambda, \mu)} & \text{for\, $2\leq p \leq \infty$}\end{array} \right.
\end{align*}
Moreover, ${\G}^d_{p,q}$, the closed subspace of $L_p(\Gamma_q)$ spanned by $$\{g_{p,q,\und{i}} : \und{i}\in \I^d\}$$
is completely complemented with the constant depending only on $d$ and $q$. 
\end{thm}
\begin{proof}
We consider the map $u_{p,q}$ as before with a suitable modification.
For $1\leq p\leq 2$
\begin{align*}
u_{p,q} : \; & \K^d_p(\lambda, \mu) \rightarrow L_p(\M \bar{\otimes} \Gamma_q)\\
& \;\;\;\; x \;\; \mapsto \sum_{\und{i}\in \I^d_m}x_{\und{i}}\otimes g_{p,q,\und{i}},
\end{align*}
and for $2\leq p\leq \infty$
\begin{align*}
u_{p,q} : \; & \J^d_p(\lambda, \mu) \rightarrow L_p(\M \bar{\otimes} \Gamma_q)\\
& (I_{L_p(\M)}\otimes P_d)x \;\; \mapsto \sum_{\und{i}\in \I^d_m}x_{\und{i}}\otimes g_{p,q,\und{i}}.
\end{align*}
We observe that
	$$\sum_{\und{i}\in \I^d_m}x_{\und{i}}\otimes g_{p,q,\und{i}} = (I_{L_p(\M)} \otimes D^{\frac{1}{2p}})(I_{L_p(\M)}\otimes W)(\xi)
	(I_{L_p(\M)} \otimes D^{\frac{1}{2p}}),$$
where $$\xi = \sum_{\und{i}\in \I^d_m} x_{\und{i}} \otimes \lambda_{\und{i}}e_{\und{i}}.$$
Thus, for $2\leq p\leq \infty$ we have
$$u_{p,q}(x) = (I_{L_p(\M)} \otimes D^{\frac{1}{2p}})(I_{L_p(\M)}\otimes W)[(I_{L_p(\M)}\otimes P^{-1}_d)(\xi)](I_{L_p(\M)} \otimes D^{\frac{1}{2p}}).$$

Our plan is as follows. First, we will show that $u_{\infty, q}$, $u_{2,q}$ and $u_{1,q}$ are bounded with constants depending only on $d$ and $q$.
Then by Corollary \ref{cor-interpolation1} and complex interpolation
$u_{p,q}$'s $(1\leq p \leq \infty)$ are bounded with constants depending only on $d$ and $q$,
since it is clear that $u_{p,q}$'s are essentially the same map for $1\leq p \leq 2$ and $2\leq p \leq \infty$,
respectively, in the sense of interpolation theory.
Then, we will apply a duality argument to show that $u_{p,q}$'s $(1\leq p \leq \infty)$ are actually isomorphisms
with constants depending only on $d$ and $q$.
Indeed, for $x = (x_{\und{i}})_{\und{i}\in \I^d_m}$, $z = (z_{\und{i}})_{\und{i}\in \I^d_m}$ and $1\leq p \leq 2$
	\begin{align}\label{duality3}
	\lefteqn{\left\langle u_{p,q}(x), u_{p',q}(z) \right\rangle}\\
	& = \sum_{\und{i}, \und{j}\in \I^d_m}\text{tr}_{\M}(x^*_{\und{i}}z_{\und{j}}) \text{tr}_{\Gamma_q}
	(D^{\frac{1}{2p}}W^*(\lambda_{\und{i}}e_{\und{i}})D^{\frac{1}{2p}}D^{\frac{1}{2p'}}W(P^{-1}_d\lambda_{\und{j}}e_{\und{j}})D^{\frac{1}{2p'}})\nonumber\\
	& = \sum_{\und{i},\und{j}\in \I^d_m}
	\lambda^{-\frac{1}{p}}_{\und{i}}\lambda^{-\frac{1}{p'}}_{\und{j}}\mu^{\frac{1}{p}}_{\und{i}}\mu^{\frac{1}{p'}}_{\und{j}}
	\text{tr}_{\M}(x^*_{\und{i}}z_{\und{j}}) \phi_q(W^*(\lambda_{\und{i}}e_{\und{i}})W(P^{-1}_d\lambda_{\und{j}}e_{\und{j}}))\nonumber\\
	& = \sum_{\und{i}\in \I^d_m} \lambda_{\und{i}}\mu_{\und{i}}\text{tr}_{\M}(x^*_{\und{i}}z_{\und{i}})
	= \left\langle x, z \right\rangle \nonumber
	\end{align}
and for any $Y \in L_p(\M \bar{\otimes} \Gamma_q)$ and $1\leq p\leq \infty$
	\begin{align}\label{duality4}
	\left\langle u_{p,q}(x), Y \right\rangle = \left\langle u_{p,q}(x), u_{p',q} u^*_{p,q}(Y) \right\rangle.
	\end{align}
If we know that $u_{p,q}$'s $(1\leq p \leq \infty)$ are bounded with constants depending only on $d$ and $q$,
then we can conclude that $u_{p,q}$'s $(1\leq p \leq \infty)$ are isomorphisms with constants depending only on $d$ and $q$
using \eqref{duality3} and \eqref{duality4}.

Now we consider the boundedness of $u_{\infty, q}$ following the approach of Nou (\cite{N04}).
For $x = (x_{\und{i}})_{\und{i}\in \I^d_m}$ we have by \eqref{eq-Wick} and \eqref{Delta-info}
\begin{align*}
u_{\infty, q}(x) & = \sum^d_{k=0} I_\M \otimes [U_k (I_k\otimes S)(I_k\otimes S) R^*_{d,d-k}](\tilde{\xi})\\
& = \sum^d_{k=0} I_\M \otimes [U_k (I_k\otimes S)(I_k\otimes J)(I_k\otimes (A^{-\frac{1}{2}})^{\otimes d-k}) R^*_{d,d-k}](\tilde{\xi}),
\end{align*}
where $\tilde{\xi} = (I_{\M}\otimes P^{-1}_d)(\xi)$ and $I_k$ is the formal identity on $\Hi^k$.

Since $\norm{U_k (I_k\otimes S)}_{cb}\leq C_q$ and $\norm{I_k\otimes J}_{cb} \leq 1$ we have
\begin{align*}
\lefteqn{\norm{u_{\infty,q}(x)}_{\M \otimes_{\min} \Gamma_q}}\\
& \leq C_q \sum^d_{k=0} \norm{I_\M \otimes (I_k\otimes (A^{-\frac{1}{2}})^{\otimes d-k})
R^*_{d,d-k}(\tilde{\xi})}_{\M \otimes_{\min}(H^{\otimes k}_c \otimes_h H^{\otimes d-k}_r)}.
\end{align*}
By \eqref{P-n-isometry} and \eqref{norm-P-n-inverse} we have
\begin{align*}
\lefteqn{\norm{I_\M \otimes (I_k\otimes (A^{-\frac{1}{2}})^{\otimes d-k})
R^*_{d,d-k}(\tilde{\xi})}_{\M \otimes_{\min}(H^{\otimes k}_c \otimes_h H^{\otimes d-k}_r)}}\\
& = \norm{I_\M \otimes [(P^{\frac{1}{2}}_k \otimes P^{\frac{1}{2}}_{d-k})(I_k\otimes (A^{-\frac{1}{2}})^{\otimes d-k})
R^*_{d,d-k}](\tilde{\xi})}_{\M \otimes_{\min}(\Hi^{\otimes k}_c \otimes_h \Hi^{\otimes d-k}_r)}\\
& \sim_{c_{q,d}} \norm{I_\M \otimes [(P_k \otimes P_{d-k})(I_k\otimes (A^{-\frac{1}{2}})^{\otimes d-k})
R^*_{d,d-k}](\tilde{\xi})}_{\M \otimes_{\min}\Hi^{\otimes k}_c \otimes_h \Hi^{\otimes d-k}_r}.
\end{align*}
Since $P_{d-k}$ and $(A^{-\frac{1}{2}})^{\otimes d-k}$ commute we have by \eqref{P-n-R-nk} that
\begin{align*}
(P_k \otimes P_{d-k})(I_k\otimes (A^{-\frac{1}{2}})^{\otimes d-k})R^*_{d,d-k}
& = (I_k\otimes (A^{-\frac{1}{2}})^{\otimes d-k})(P_k \otimes P_{d-k})R^*_{d,d-k}\\
& = (I_k\otimes (A^{-\frac{1}{2}})^{\otimes d-k})P_d.
\end{align*}
Thus, by \eqref{action-A} we have 
\begin{align*}
\lefteqn{\norm{I_\M \otimes (I_k\otimes (A^{-\frac{1}{2}})^{\otimes d-k})
R^*_{d,d-k}(\tilde{\xi})}_{\M \otimes_{\min}(H^{\otimes k}_c \otimes_h H^{\otimes d-k}_r)}}\\
& \sim_{c_{q,d}} \norm{[I_\M \otimes (I_k\otimes (A^{-\frac{1}{2}})^{\otimes d-k})](\xi)}_{\M
\otimes_{\min}(\Hi^{\otimes k}_c \otimes_h \Hi^{\otimes d-k}_r)}\\
& = \norm{\sum_{\und{i}\in \I^d_m} x_{\und{i}}\otimes
\lambda_{\und{i}^k}e_{\und{i}^k} \otimes \mu_{\und{i}^{d-k}}e_{\und{i}^{d-k}}}_{\M \otimes_{\min}(\Hi^{\otimes k}_c \otimes_h \Hi^{\otimes d-k}_r)}\\
& = \R\C^{d,k}_{\infty}(\lambda,\mu; x).
\end{align*}
Consequently, we get
$$\norm{u_{\infty,q}(x)}_{\M \otimes_{\min} \Gamma_q} \lesssim_{c_{q,d}} \norm{x}_{\J^d_\infty(\lambda,\mu)}.$$

When $p=2$ the calculation is straightforward. We can easily check that
\begin{align*}
\norm{u_{2,q}(x)}_{L_2(\M \bar{\otimes} \Gamma_q)} & = \norm{\sum_{\und{i}\in\I^d_m}x_{\und{i}}\otimes
\lambda^{\frac{1}{2}}_{\und{i}}\mu^{\frac{1}{2}}_{\und{i}}e_{\und{i}}}_{L_2(\M)\otimes_2 H^{\otimes d}}\\
& = \norm{(I_{L_2(\M)} \otimes P^{\frac{1}{2}}_d)\sum_{\und{i}\in\I^d_m}x_{\und{i}}\otimes
\lambda^{\frac{1}{2}}_{\und{i}}\mu^{\frac{1}{2}}_{\und{i}}e_{\und{i}}}_{L_2(\M)\otimes_2 \Hi^{\otimes d}}\\
& = \norm{(I_{L_2(\M)} \otimes P^{\frac{1}{2}}_d)x}_{\J^d_2(\lambda,\mu)}.
\end{align*}
Since $$\norm{(I_{L_2(\M)} \otimes P^{\frac{1}{2}}_d)x}_{\J^d_2(\lambda,\mu)} \sim_{c_{q,d}} \norm{(I_{L_2(\M)} \otimes P_d)x}_{\J^d_2(\lambda,\mu)}
\sim_{c_{q,d}} \norm{x}_{\K^d_2(\lambda,\mu)}$$ by \eqref{norm-P-n} and \eqref{norm-P-n-inverse} we get our desired estimate.

The boundedness of $u_{1,q}$ can be checked by the identical argument as in the proof of Theorem \ref{thm-free}, so we are done here.

\end{proof}

\section{The case of CAR and CCR generators}\label{sec-CAR}

We will simply write $a_{k}$, $a_{p,k}$ and $\A^d_p$ instead of $g_{\pm 1,k}$, $g_{p,\pm1 ,k}$ and $\G^d_{p,\pm 1}$, respectively.
For any $\und{i} \in \n^d_m$ we define $a_{p, \und{i}}$ similarly. 
We start this section with the operator space structure of $\A^1_p$.

\begin{thm}\label{degree-1}
Let $m\in \n$ and $x = (x_k)_{1\leq k \leq m} \subseteq L_p(\M).$ Then we have
\begin{align*}
\norm{\sum_{1\leq k \leq m}x_k \otimes a_{p,k}}_{L_p(\M \bar{\otimes}\Gamma_{\pm 1})}
& \sim_{c_p} \left\{\begin{array}{ll}\norm{x}_{\K^1_p(\lambda, \mu)} & \text{for\, $1\leq p \leq 2$}\\
\norm{x}_{\J^1_p(\lambda, \mu)} & \text{for\, $2\leq p < \infty$.}\end{array} \right.
\end{align*}
Moreover, $L_p(\M; \A_p)$ is complemented in $L_p(\M \bar{\otimes}\Gamma_{\pm 1})$.

In particular, we have $\A^1_p \cong RC_p(\lambda, \mu)$ ($1\leq p<\infty$) completely isomorphically,
and $\A^1_p$ ($1< p <\infty$) is completely complemented in $L_p(\Gamma_{\pm 1})$. All constants here depend only on $p$.
\end{thm}
\begin{proof}
When $2\leq p <\infty$ we apply noncommutative Burkholder inequality (also see Theorem 4.1 of \cite{X-Rp-Embedding}).
For $1<p<2$ we follow the proof of Theorem E in \cite{JPX-FreeChaos} to obtain an upper bound.
The lower bound follows by the duality and the previous result.

For the case $p=1$ we can get the result with a slight modification of the proof of Theorem 7.1 in \cite{J-Araki}.
Let $\psi$, $u_n$ and $N$ be as in section \ref{subsec-matrix}. Let $N_n = M_{2^n}\otimes N^{\otimes n}$ and
$D_\psi$, $D_n$ be the densities of $\psi$ and $\tau_n \otimes (\frac{\psi}{4m})^{\otimes n}$, respectively.
Then by \eqref{translation} and \eqref{matrix-model-inclusion} we have
\begin{align*}
\lefteqn{\norm{\sum^m_{k = 1}x_k \otimes a_{1,k}}_{L_1(\M \bar{\otimes}\Gamma_{\pm 1})}}\\
& =  \lim_n \norm{\sum^m_{k = 1} \frac{\sqrt{\lambda^2_k + \mu^2_k}}{2} x_k \otimes
D^{\frac{1}{2}}_n u_n\big((\delta_k + i\delta_{-k})\otimes e_{12}\big) D^{\frac{1}{2}}_n}_{L_1(\M \bar{\otimes}N_n)}.
\end{align*}
The only difference now is that we are using $(\delta_k + i\delta_{-k})\otimes e_{12}$
instead of $\delta_k \otimes e_{12}$ in (7.3) of \cite{J-Araki}. Since
\begin{align*}
\lefteqn{\text{tr}_N ( [D^{\frac{1}{2}}_{\psi}\big((\delta_l - i\delta_{-l})\otimes e_{21}\big)]^*
D^{\frac{1}{2}}_{\psi} \big((\delta_k + i\delta_{-k})\otimes e_{12}\big))}\\ & =
\text{tr}_N (\big((\delta_l - i\delta_{-l})\otimes e_{21}\big) D_{\psi} \big((\delta_k + i\delta_{-k})\otimes e_{12}\big) )\\
& = 2\delta_{kl}(2-\sigma_k) = \delta_{kl}\frac{4\lambda^2_{k}}{\lambda^2_{k} + \mu^2_{k}}.
\end{align*}
by the same argument we get
\begin{align*}
\lefteqn{\norm{\sum^m_{k = 1}x_k \otimes a_{1,k}}_{L_1(\M \bar{\otimes}\Gamma_{\pm 1})}}\\
& \sim \inf_{x_k = c_k + d_k} \Bigg\{ \norm{\Big(\sum^m_{k = 1} \lambda^2_{k}c^*_k c_k \Big)^{\frac{1}{2}}}_{L_1(\M)}
+ \norm{\Big(\sum^m_{k = 1} \mu^2_{k}d_k d^*_k \Big)^{\frac{1}{2}}}_{L_1(\M)}\Bigg\}\\
& = \norm{(x_k)^m_{k=1}}_{\K^1_1(\lambda, \mu)}.
\end{align*}

\end{proof}

For $x = (x_{\und{i}})_{\und{i}\in \n^d_m} \subseteq L_p(\M)$ and $\rb = (\rb_1, \cdots, \rb_d) \in \{c, r\}^d$ we denote
$$\R\C_p^{\rb}(\lambda, \mu; x) = \norm{\sum_{\und{i}\in \n^d_m}x_{\und{i}}\otimes
h_{p, \rb_1, i_1}\otimes \cdots \otimes h_{p, \rb_d, i_d}}_{L_p(\M ; F^{\rb}_p)},$$ where
$$h_{p, c, i} = \lambda^{\frac{1}{p'}}_i \mu^{\frac{1}{p}}_i e_{i1} \in C_p,\;\,
h_{p, r, i} = \lambda^{\frac{1}{p}}_i \mu^{\frac{1}{p'}}_i e_{1i}\in R_p,$$
$$F^{\rb}_p = F_1 \otimes_h \cdots \otimes_h F_d \subseteq S_p(\ell^{\otimes d}_2)$$
and $$F_k = \left\{ \begin{array}{ll} C_p & \text{if}\;\, \rb_k = c \\ R_p & \text{if}\;\, \rb_k = r. \end{array} \right.$$
We define the corresponding symmetric $\K$- and $J$-functional spaces
$\s\K^d_p(\lambda, \mu)$ $(1\leq p\leq 2)$ and $\s\J^d_p(\lambda, \mu)$ $(2\leq p\leq \infty)$
as the closures of finite tuples in $L_p(\M)$ indexed by $\n^d$ with respect to the following norms. 
$$\norm{x}_{\s\K^d_p(\lambda,\mu)}=\inf_{\rb \in \{c,r\}^d} \R\C_p^{\rb}(\lambda, \mu; x^\rb),$$
where the infimum runs over all possible $x_{\und{i}} = \sum_{\rb \in \{c,r\}^d}x^{\rb}_{\und{i}}$ and
$x^\rb = (x^\rb_{\und{i}})_{\und{i}\in \n^d_m}$ for $\rb \in \{c,r\}^d$, and
$$\norm{x}_{\s\J^d_p(\lambda, \mu)} = \max_{\rb \in \{c,r\}^d} \R\C_p^{\rb}(\lambda, \mu; x).$$

The following is the tensor product extension of Theorem \ref{degree-1}.
\begin{cor}\label{lem-tensor}
For $x = (x_{\und{i}})_{\und{i}\in \n^d_m} \subseteq L_p(\M)$ and $1\leq p < \infty$ we have the following equivalence.
\begin{align*}
\norm{\sum_{\und{i}\in \n^d_m}x_{\und{i}}\otimes a_{p,i_1}\otimes \cdots \otimes
a_{p,i_d}}_{L_p(\M \bar{\otimes}\Gamma_{\pm 1}^{\otimes d})} \sim_{c_{p,d}}
\left\{\begin{array}{ll}\norm{x}_{\s\K^d_p(\lambda, \mu)} & \text{for\, $1\leq p \leq 2$}\\
\norm{x}_{\s\J^d_p(\lambda, \mu)} & \text{for\, $2\leq p < \infty$.}\end{array} \right.
\end{align*}
\end{cor}

\begin{proof}
We assume that $1\leq p\leq 2$ and $d=2$. The proof for the other cases are the same.
By Theorem \ref{degree-1} we have
\begin{align*}
\lefteqn{\norm{\sum^m_{i, j =1}x_{ij} \otimes a_{p,i}\otimes a_{p,j}}_{L_p(\M \bar{\otimes}\Gamma_{\pm 1} \bar{\otimes}\Gamma_{\pm 1})}}\\
& \sim_{c_p} & \inf_{x_{ij} = c_{ij} + d_{ij}} \Bigg\{ \norm{\Big(\sum^m_{j =1}(\sum^m_{i=1} c_{ij} \otimes a_{p,i}) \otimes \lambda^{\frac{1}{p}}_{j}\mu^{\frac{1}{p'}}_j e_{j1} \Big)^{\frac{1}{2}}}_{L_p(\M \bar{\otimes}\Gamma_{\pm 1} \bar{\otimes}B(\ell_2))}\\
& & + \norm{\Big(\sum^m_{j =1}(\sum^m_{i=1} d_{ij} \otimes a_{p,i})
\otimes \lambda^{\frac{1}{p'}}_{j}\mu^{\frac{1}{p}}_j e_{1j} \Big)^{\frac{1}{2}}}_{L_p(\M \bar{\otimes}\Gamma_{\pm 1} \bar{\otimes}B(\ell_2))}\Bigg\}.
\end{align*}
If we apply Theorem \ref{degree-1} for $\M \bar{\otimes}B(\ell_2)$, then we get the desired result.
\end{proof}

\subsection{Polynomials of degree 2}
In this section we first focus on the case $d=2$ to get the idea for polynomials of higher degree.

\begin{thm}\label{thm-main-degree-2}
Let $x = (x_{ij})^m_{i,j=1} \subseteq L_p(\M)$ satisfying
\begin{equation}\label{symmetry-degree-2}
x_{ij} = -x_{ji}\;\, \text{for}\;\, q=-1
\end{equation}
and $$x_{ij} = x_{ji} \;\, \text{for}\;\, q=1$$
for all $1\leq i, j \leq m$. Then we have the following equivalence.
\begin{align*}
\norm{\sum^m_{i,j=1}x_{ij}\otimes D^{\frac{1}{2p}}_{\pm 1}a_i a_j D^{\frac{1}{2p}}_{\pm 1}}_{L_p(\M \bar{\otimes}\Gamma_{\pm 1})}
& \sim_{c_{p}} \left\{\begin{array}{ll}\norm{x}_{\s\K^2_p(\lambda, \mu)} & \text{for\, $1\leq p \leq 2$}\\
\norm{x}_{\s\J^2_p(\lambda, \mu)} & \text{for\, $2\leq p < \infty$.}\end{array} \right.
\end{align*}
\end{thm}

\begin{proof}
Let $u_n$ and $v_k$ be as in \eqref{def-u_n} and \eqref{def-v_k}, respectively. By \eqref{matrix-model-inclusion} we have
\begin{align*}
\lefteqn{\norm{\sum^m_{i,j=1}x_{ij}\otimes D^{\frac{1}{2p}}_{\pm 1}a_i a_j
D^{\frac{1}{2p}}_{\pm 1}}_{L_p(\M \bar{\otimes}\Gamma_{\pm 1})}}\\
& = \lim_n \norm{\sum^m_{i,j=1}x_{ij}\otimes D^{\frac{1}{2p}}_n (u_n(b_i) u_n(b_j)) D^{\frac{1}{2p}}_n}_{L_p(\M \bar{\otimes}N_n)},
\end{align*}
where $b_k = \frac{\sqrt{\lambda^2_k + \mu^2_k}}{2}(\delta_k + i\delta_{-k})\otimes e_{12}$.

Note that we have $$u_n(b_i)u_n(b_j) = \frac{4m}{n}\sum^n_{k,l=1}v_k v_l \otimes \pi_k(b_i)\pi_l(b_j)$$ and
$\psi_n(v_k v_l \otimes \pi_k(b_i)\pi_l(b_j)) = 0$ for $k\neq l$. Thus, we have
\begin{align*}
\lefteqn{\norm{\sum^m_{i,j=1}x_{ij}\otimes D^{\frac{1}{2p}}_n (u_n(b_i) u_n(b_j))D^{\frac{1}{2p}}_n}_p}\\
& = \frac{4m}{n}\norm{\sum^m_{i,j=1}x_{ij}\otimes D^{\frac{1}{2p}}_n \Big(\sum^n_{k,l=1}v_k v_l \otimes \pi_k(b_i)\pi_l(b_j)\Big)
D^{\frac{1}{2p}}_n}_p.
\end{align*}
We set $$Q_n = \sum_{\substack{1\leq i, j \leq m\\ 1\leq k, l \leq n}}x_{ij}\otimes D^{\frac{1}{2p}}_n
\Big(v_k v_l \otimes\pi_k(b_i)\pi_l(b_j)\Big) D^{\frac{1}{2p}}_n.$$

First, we show the lower bound by an averaging technique. Let us fix a subset $A \subseteq \{1,\cdots,n\}$. We introduce
$$v_k(\vep) = \left\{ \begin{array}{ll}\vep_1 v_k & \text{for}\;\, k\in A\\ \vep_2 v_k & \text{for}\;\, k\notin A \end{array}\right.,$$
where $\vep_1$ and $\vep_2$ are two independent Bernoulli variables with $\text{Prob}(\vep_i = \pm 1) = \frac{1}{2}.$
Then
\begin{align*}
& \int^1_0 \int^1_0 \vep_1(t_1)\vep_2(t_2) \sum_{\substack{1\leq i, j \leq m\\ 1\leq k, l \leq n}}x_{ij}\otimes
v_k(\vep) v_l(\vep) \otimes\pi_k(b_i)\pi_l(b_j)dt_1 dt_2\\
& = \sum_{\substack{1\leq i, j \leq m\\k\in A,\; l\notin A}}x_{ij}\otimes v_k v_l \otimes\pi_k(b_i)\pi_l(b_j)
+ \sum_{\substack{1\leq i, j \leq m\\k\notin A,\; l\in A}}x_{ij}\otimes v_k v_l \otimes\pi_k(b_i)\pi_l(b_j)\\
& = 2 \sum_{\substack{1\leq i, j \leq m\\k\in A,\; l\notin A}}x_{ij}\otimes v_k v_l \otimes\pi_k(b_i)\pi_l(b_j)
\end{align*}
since (when $q=-1$)
\begin{align*}
\sum_{\substack{1\leq i, j \leq m\\k\in A,\; l\notin A}}x_{ij}\otimes v_k v_l \otimes\pi_k(b_i)\pi_l(b_j)
& = \sum_{\substack{1\leq i, j \leq m\\k\notin A,\; l\in A}}x_{ij}\otimes v_l v_k \otimes\pi_l(b_i)\pi_k(b_j)\\
& = -\sum_{\substack{1\leq i, j \leq m\\k\notin A,\; l\in A}}x_{ij}\otimes v_k v_l \otimes\pi_k(b_j)\pi_l(b_i)\\
& = \sum_{\substack{1\leq i, j \leq m\\k\notin A,\; l\in A}}x_{ij}\otimes v_k v_l \otimes\pi_k(b_i)\pi_l(b_j).
\end{align*}
The case $q=1$ is the same.

Clearly the tuple $\big( (v_k(\vep))_{k\in A}, (v_k(\vep))_{k\notin A}\big)$ has the same $*$-distribution as
the tuple $\big( (v_k)_{k\in A}, (v_k)_{k\notin A}\big)$.
Thus, by combining the above observations we get
\begin{equation}\label{lower-bound-deg-2}
\norm{\sum_{\substack{1\leq i, j \leq m\\k\in A,\; l\notin A}}x_{ij}\otimes
D^{\frac{1}{2p}}_n\Big( v_k v_l \otimes \pi_k(b_i)\pi_l(b_j)\Big) D^{\frac{1}{2p}}_n}_p \leq \frac{1}{2}\norm{Q_n}_p.
\end{equation}

Now we assume that $n$ is even and consider a specific choice of $A$, namely $$A = \{1,\cdots, \frac{n}{2}\}.$$
Then $k\in A,\; l \notin A$ means $1\leq k \leq \frac{n}{2} < l \leq n$.
In this case we decompose $v_k$ and $v_j$ as follows. $$v_k = w_k\otimes 1, \;\, v_l = u \otimes w_l,$$
where $w_k, w_l$ and $u \in M_{2^{n/2}}$. Moreover, $\pi_k(b_i)\pi_l(b_j)$ can be understood as $\pi_k(b_i)\otimes \pi_l(b_j)$
since $\pi_k(b_i) \in N^{\otimes \frac{n}{2}} \otimes 1$ and $\pi_l(b_j) \in 1 \otimes N^{\otimes \frac{n}{2}}$.
Since $D_n = I\otimes D^{\otimes n}_{\psi/4m}$ we can write $D_n = D_{n/2} \otimes D_{n/2}$. Now we have 
\begin{align*}
\lefteqn{\sum^m_{i, j=1}\sum_{1\leq k \leq \frac{n}{2} < l \leq n}x_{ij}\otimes
D^{\frac{1}{2p}}_n\Big[v_k v_l \otimes\pi_k(b_i)\pi_l(b_j)\Big]D^{\frac{1}{2p}}_n}\\
& = \sum^m_{i, j=1}\sum_{1\leq k \leq \frac{n}{2} < l \leq n}x_{ij}\otimes
D^{\frac{1}{2p}}_n \Big[(w_k \otimes 1) (u \otimes w_l) \otimes \pi_k(b_i)\pi_l(b_j)\Big]D^{\frac{1}{2p}}_n\\
& = \sum^m_{i, j=1}x_{ij}\otimes D^{\frac{1}{2p}}_{n/2}\Big[\sum^{n/2}_{k=1}w_k u\otimes\pi_k(b_i)\Big]D^{\frac{1}{2p}}_{n/2}
\otimes D^{\frac{1}{2p}}_{n/2}\Big[\sum^n_{l=n/2+1}w_l\otimes \pi_l(b_j)\Big]D^{\frac{1}{2p}}_{n/2}.
\end{align*}
Since $w_k$ and $u$ commute $(w_k u)_{1\leq k \leq \frac{n}{2}}$ and $(w_l)_{\frac{n}{2}< l \leq n}$ satisfy the relations in \eqref{relation-v_k}. 
Thus, $$\sqrt{\frac{4m}{n/2}}\sum^{n/2}_{k=1}w_k u\otimes\pi_k(b_i)\;\, \text{and}\;\, \sqrt{\frac{4m}{n/2}}\sum^n_{l=n/2+1}w_l\otimes \pi_l(b_j)$$
converges in $*$-distribution to $a_i$ and $a_j$, respectively, and consequently we have
\begin{align}\label{decoupling-deg-2}
\lefteqn{\norm{\sum^m_{i,j=1}x_{ij}\otimes a_{i,p} \otimes a_{j,p}}_p}\nonumber \\ 
& = \lim_n \frac{4m}{n/2}\norm{\sum^m_{i, j=1}\sum_{1\leq k \leq \frac{n}{2} < l \leq n}x_{ij}\otimes
D^{\frac{1}{2p}}_n\Big(v_k v_l \otimes\pi_k(b_i)\pi_l(b_j)\Big)D^{\frac{1}{2p}}_n}_p\\
& \leq 2\lim_n \frac{4m}{n}\norm{Q_n}_p = 2 \norm{\sum^m_{i,j=1}x_{ij}\otimes D^{\frac{1}{2p}}_{\pm 1}a_i a_j D^{\frac{1}{2p}}_{\pm 1}}_p. \nonumber
\end{align}
The conclusion follows by Corollary \ref{lem-tensor}.

For the upper bound we consider the discrete probability space $\Om$ of $$A\subseteq \{1,\cdots, n\}\;\, \text{with}\;\, \abs{A} = \frac{n}{2}.$$
Then for fixed $1\leq k\neq l \leq n$ we have
\begin{align*}
\text{Prob}(k\in A, l\notin A)
& = \frac{\#\{A\subseteq \{1,\cdots, n\} : \abs{A} = \frac{n}{2}, k\in A, l\notin A\}}{\#\{A\subseteq \{1,\cdots, n\} : \abs{A} = \frac{n}{2}\}}\\
& = \frac{\binom{n-2}{n/2-1}}{\binom{n}{n/2}} = \frac{n}{4(n-1)}.
\end{align*}
Thus, by the symmetry we have
\begin{align*}
\norm{Q_n}_p & = \frac{4(n-1)}{n}\norm{\sum^n_{k,l=1}\E_A (1_{k\in A, l\notin A})\sum^m_{i,j=1}
x_{ij}\otimes D^{\frac{1}{2p}}_n \Big(v_k v_l \otimes\pi_k(b_i)\pi_l(b_j)\Big) D^{\frac{1}{2p}}_n}_p\\
& \leq 4 \E_A \norm{\sum_{k\in A, l\notin A}\sum^m_{i,j=1}
x_{ij}\otimes D^{\frac{1}{2p}}_n \Big(v_k v_l \otimes\pi_k(b_i)\pi_l(b_j)\Big) D^{\frac{1}{2p}}_n}_p\\
& = 4 \norm{\sum_{1 \leq k \leq n/2 < l \leq n}\sum^m_{i,j=1}
x_{ij}\otimes D^{\frac{1}{2p}}_n \Big(v_k v_l \otimes\pi_k(b_i)\pi_l(b_j)\Big) D^{\frac{1}{2p}}_n}_p,
\end{align*}
where $\E_A$ implies the expectation with respect to $\Om$. Consequently, by \eqref{decoupling-deg-2} we get
$$\norm{\sum^m_{i,j=1}x_{ij}\otimes D^{\frac{1}{2p}}_{\pm 1}a_i a_j D^{\frac{1}{2p}}_{\pm 1}}_p
\leq 8\norm{\sum^m_{i,j=1}x_{ij}\otimes a_{p,i} \otimes a_{p,j}}_p,$$
and the conclusion follows by Corollary \ref{lem-tensor} again.
\end{proof}

\subsection{Polynomials of higher degree}
Now we consider the general case $d\geq 2$. For a permutation $\rho \in S_d$ and $\und{i} \in \n^d_m$ we denote
$(i_{\rho(1)},\cdots ,i_{\rho(d)}) \in \n^d_m$ by $\und{i}(\rho)$.

\begin{thm}\label{thm-main-degree-d}
Let $x = (x_{\und{i}})_{\und{i}\in \n^d_m} \subseteq L_p(\M)$ satisfying
\begin{equation}\label{symmetry-degree-d}
x_{\und{i}} = \text{sgn}(\rho)x_{\und{i}(\rho)}\;\, \text{for}\;\, q=-1
\end{equation}
$$x_{\und{i}} = x_{\und{i}(\rho)} \;\, \text{for}\;\, q=1$$
for any $\rho \in S_d$ and $\und{i} \in \n^d_m$. Then we have the following equivalence.
\begin{align*}
\norm{\sum_{\und{i}\in \n^d_m}x_{\und{i}}\otimes a_{p,\und{i}}}_{L_p(\M \bar{\otimes}\Gamma_{\pm 1})}
& \sim_{c_{p,d}} \left\{\begin{array}{ll}\norm{x}_{\s\K^d_p(\lambda, \mu)} & \text{for\, $1\leq p \leq 2$}\\
\norm{x}_{\s\J^d_p(\lambda, \mu)} & \text{for\, $2\leq p < \infty$.}\end{array} \right.
\end{align*}
\end{thm}

\begin{proof}
As in the proof of Theorem \ref{thm-main-degree-2} we have
\begin{align*}
\lefteqn{\norm{\sum_{\und{i}\in \n^d_m}x_{\und{i}}\otimes a_{p,\und{i}}}_{L_p(\M \bar{\otimes}\Gamma_{\pm 1})}
 = \norm{\sum_{\und{i}\in \n^d_m}x_{\und{i}}\otimes D^{\frac{1}{2p}}_{\pm 1}a_{\und{i}}D^{\frac{1}{2p}}_{\pm 1}}_p}\\
& = \lim_n \norm{\sum_{\und{i}\in \n^d_m}x_{\und{i}}\otimes D^{\frac{1}{2p}}_n
(u_n(b_{i_1})\cdots u_n(b_{i_d})) D^{\frac{1}{2p}}_n}_{L_p(\M \bar{\otimes}N_n)}\\
& = \lim_n \Big(\frac{4m}{n}\Big)^{\frac{d}{2}}\norm{\sum_{\und{i}\in \n^d_m}x_{\und{i}}\otimes 
\sum_{\und{j}\in \n^d_n} v_{\und{j}} \otimes D^{\frac{1}{2p}}_n (\pi_{j_1}(b_{i_1})\cdots \pi_{j_d}(b_{i_d})) D^{\frac{1}{2p}}_n}_p,
\end{align*}
and we set $$Q_n = \sum_{\und{i}\in \n^d_m, \und{j}\in \n^d_n} x_{\und{i}}\otimes v_{\und{j}}
\otimes D^{\frac{1}{2p}}_n (\pi_{j_1}(b_{i_1})\cdots \pi_{j_d}(b_{i_d})) D^{\frac{1}{2p}}_n$$

First, we consider the lower bound. Let $(A_1, \cdots, A_d)$ be a partition of $\{1,\cdots,n\}$.
For $1\leq k \leq d$ we introduce $$v_j(\varepsilon) = \varepsilon_k v_j\;\, \text{for all}\;\, j \in A_k,$$
where $(\varepsilon_k)^d_{k=1}$ is a collection of independent Bernoulli variables with $\text{Prob}(\vep_i = \pm 1) = \frac{1}{2}.$
Now we observe that
\begin{align*}
\lefteqn{\int_{[0,1]^d} \varepsilon_1 \cdots \varepsilon_d \sum_{\und{i}\in \n^d_m}x_{\und{i}}\otimes 
\sum_{\und{j}\in \n^d_n} v_{\und{j}}(\varepsilon) \otimes D^{\frac{1}{2p}}_n (\pi_{j_1}(b_{i_1})\cdots \pi_{j_d}(b_{i_d}))
D^{\frac{1}{2p}}_n}\\
& = \sum_{\rho \in S_d}\sum_{j_{\rho(1)}\in A_1, \cdots, j_{\rho(d)}\in A_d} \sum_{\und{i}\in \n^d_m}x_{\und{i}}\otimes v_{\und{j}} \otimes D^{\frac{1}{2p}}_n (\pi_{j_1}(b_{i_1})\cdots \pi_{j_d}(b_{i_d})) D^{\frac{1}{2p}}_n
\end{align*}
since the integration kills all the terms where all $j_k$'s belong to different partitions. 

On the other hand for a fixed $\rho \in S_d$ we have (when $q=-1$)
\begin{align*}
\lefteqn{\sum_{j_{\rho(1)}\in A_1, \cdots, j_{\rho(d)}\in A_d} \sum_{\und{i}\in \n^d_m}x_{\und{i}}\otimes v_{\und{j}} \otimes 
\pi_{j_1}(b_{i_1})\cdots \pi_{j_d}(b_{i_d})}\\
& = \sum_{j_1\in A_1, \cdots, j_d\in A_d} \sum_{\und{i}\in \n^d_m}x_{\und{i}}\otimes v_{\und{j}(\rho^{-1})} \otimes  \pi_{j_{\rho^{-1}(1)}}(b_{i_1})\cdots \pi_{j_{\rho^{-1}(d)}}(b_{i_d})\\
& = \sum_{j_1\in A_1, \cdots, j_d\in A_d} \text{sgn}(\rho)\sum_{\und{i}\in \n^d_m}x_{\und{i}}\otimes v_{\und{j}} \otimes  \pi_{j_1}(b_{i_{\rho^{-1}(1)}})\cdots \pi_{j_d}(b_{i_{\rho^{-1}(d)}})\\
& = \sum_{j_1\in A_1, \cdots, j_d\in A_d}\sum_{\und{i}\in \n^d_m} \text{sgn}(\rho)x_{\und{i}(\rho)}\otimes v_{\und{j}} \otimes  \pi_{j_1}(b_{i_1})\cdots \pi_{j_d}(b_{i_d})\\
& = \sum_{j_1\in A_1, \cdots, j_d\in A_d}\sum_{\und{i}\in \n^d_m} x_{\und{i}}\otimes v_{\und{j}} \otimes \pi_{j_1}(b_{i_1})\cdots \pi_{j_d}(b_{i_d}).
\end{align*}
The case $q=1$ is similar.

Since $(v_j(\varepsilon))^n_{j=1}$ has the same $*$-distribution with $(v_j)^n_{j=1}$ by combining the above two observations we get
$$d!\norm{\sum_{\und{i}\in \n^d_m}\sum_{j_1\in A_1, \cdots, j_d\in A_d} x_{\und{i}}\otimes v_{\und{j}} \otimes
D^{\frac{1}{2p}}_n(\pi_{j_1}(b_{i_1})\cdots \pi_{j_d}(b_{i_d}))D^{\frac{1}{2p}}_n}_p \leq \norm{Q_n}_p.$$

Now we assume that $n$ is a multiple of $d$ and choose
$$A_k = \Big\{\frac{(k-1)n}{d}+1, \cdots, \frac{kn}{d}\Big\}\;\, \text{for}\;\, 1\leq k\leq d.$$
As in the case $d=2$ we can write $v_j = u^{\otimes k-1}\otimes w_{j,k} \otimes 1$ for $j\in A_k$ and $u, w_{j,k} \in M_{2^{n/d}}$.
Then $$v_{\und{j}} = w_{j_1, 1}u^{d-1} \otimes v_{j_2, 2}u^{d-2} \otimes \cdots \otimes v_{j_d,d}.$$
Moreover, we have $\pi_{j_1}(b_{i_1})\cdots \pi_{j_d}(b_{i_d}) = \pi_{j_1}(b_{i_1})\otimes \cdots \otimes \pi_{j_d}(b_{i_d})$
for $j_1\in A_1, \cdots, j_d\in A_d$ and $D_n = D_{n/d}\otimes \cdots \otimes D_{n/d}$, so that
\begin{align*}
\lefteqn{\sum_{\und{i}\in \n^d_m}\sum_{j_1\in A_1, \cdots, j_d\in A_d} x_{\und{i}}\otimes v_{\und{j}} \otimes
D^{\frac{1}{2p}}_n(\pi_{j_1}(b_{i_1})\cdots \pi_{j_d}(b_{i_d}))D^{\frac{1}{2p}}_n}\\
& = \sum_{\und{i}\in \n^d_m}x_{\und{i}}\otimes \bigotimes^d_{k=1} D^{\frac{1}{2p}}_{n/d}\Big[\sum_{j_k \in A_k} v_{j_k, k}u^{d-k}\otimes
\pi_{j_k}(b_{i_k})\Big]D^{\frac{1}{2p}}_{n/d}
\end{align*}

Since $v_{j_k, k}$ and $u$ commute $(v_{j_k, k}u^{d-k})_{j_k \in A_k}$'s satisfy the relations in \eqref{relation-v_k}
$$\sqrt{\frac{4m}{n/d}}\sum_{j_k \in A_k} v_{j_k, k}u^{d-k}\otimes \pi_{j_k}(b_{i_k})$$ converges in $*$-distribution to
$a_{i_k}$ by the central limit procedure. Thus, we have
\begin{align*}
\lefteqn{\norm{\sum_{\und{i}\in \n^d_m}x_{\und{i}}\otimes a_{p,i_1} \otimes \cdots
\otimes a_{p,i_d}}_{L_p(\M \bar{\otimes} \Gamma^{\otimes d}_{\pm1})}}\\
& = \lim_n \Big(\frac{4m}{n/d}\Big)^{\frac{d}{2}}\norm{\sum_{\und{i}\in \n^d_m}x_{\und{i}}\otimes 
\bigotimes^d_{k=1}D^{\frac{1}{2p}}_{n/d}\Big[\sum_{j_k \in A_k} v_{j_k, k}u^{d-k}\otimes 
\pi_{j_k}(b_{i_k})\Big]D^{\frac{1}{2p}}_{n/d} }_{L_p(\M \bar{\otimes}N_n)}\\
& = d^{\frac{d}{2}}\lim_n \Big(\frac{4m}{n}\Big)^{\frac{d}{2}}\norm{\sum_{\und{i}\in \n^d_m}\sum_{j_1\in A_1, \cdots, j_d\in A_d}
x_{\und{i}}\otimes v_{\und{j}} \otimes D^{\frac{1}{2p}}_n(\pi_{j_1}(b_{i_1})\cdots \pi_{j_d}(b_{i_d}))D^{\frac{1}{2p}}_n}_p\\
& \leq \frac{d^{\frac{d}{2}}}{d!}\norm{\sum_{\und{i}\in \n^d_m}x_{\und{i}}\otimes a_{p,\und{i}}}_{L_p(\M \bar{\otimes}\Gamma_{\pm 1})}.
\end{align*}
Then by Corollary \ref{lem-tensor} we get the lower bound.

For the upper bound we consider a fixed $(j_1, \cdots, j_d) \in \n^d_m$
and the discrete probability space $\Om$ of partitions $\widetilde{A} = (A_1, \cdots, A_d)$ of $\{1,\cdots,n\}$ with $\abs{A_k} = \frac{n}{d}$.
Then we have 
\begin{align*}
\text{Prob}(j_1 \in A_1,\cdots, j_d \in A_d) & = \frac{\#\{ (A'_1, \cdots, A'_d) : A'_1 \cup \cdots \cup A'_d = \{1,\cdots,n-d\} \}}
{\#\{ (A_1, \cdots, A_d) : A_1 \cup \cdots \cup A_d = \{1,\cdots,n\} \}}\\
& = \frac{\binom{n-d}{n/d-1, \cdots, n/d-1}}{\binom{n}{n/d, \cdots, n/d}} = \frac{n^d(n-d)!}{d^d n!}.
\end{align*}
This implies
\begin{align*}
\lefteqn{\norm{Q_n}_p}\\ & =  \frac{d^d n!}{n^d(n-d)!}\norm{\sum_{\substack{\und{j}\in \n^d_n\\ \und{i}\in \n^d_m}}
\E_{\widetilde{A}}( 1_{j_1\in A_1,\cdots, j_d\in A_d}) x_{\und{i}}\otimes v_{\und{j}} \otimes
D^{\frac{1}{2p}}_n (\pi_{j_1}(b_{i_1})\cdots \pi_{j_d}(b_{i_d})) D^{\frac{1}{2p}}_n}_p\\
& \leq d^d \E_{\widetilde{A}}\norm{\sum_{j_1\in A_1,\cdots, j_d\in A_d}
\sum_{\und{i}\in \n^d_m}x_{\und{i}}\otimes v_{\und{j}} \otimes D^{\frac{1}{2p}}_n (\pi_{j_1}(b_{i_1})\cdots \pi_{j_d}(b_{i_d})) D^{\frac{1}{2p}}_n}_p,
\end{align*}
where $\E_{\widetilde{A}}$ is the expectation with respect to $\Om$, and consequently
\begin{align*}
\lefteqn{\norm{\sum_{\und{i}\in \n^d_m}x_{\und{i}}\otimes a_{p,\und{i}}}_{L_p(\M \bar{\otimes}\Gamma_{\pm 1})}
= \lim_n \Big(\frac{4m}{n}\Big)^{\frac{d}{2}}\norm{Q_n}_p}\\
& \leq d^d \lim_n \Big(\frac{4m}{n}\Big)^{\frac{d}{2}}\norm{\sum_{\und{i}\in \n^d_m}\sum_{j_1\in A_1, \cdots, j_d\in A_d}
x_{\und{i}}\otimes v_{\und{j}} \otimes D^{\frac{1}{2p}}_n(\pi_{j_1}(b_{i_1})\cdots \pi_{j_d}(b_{i_d}))D^{\frac{1}{2p}}_n}_p\\
& = d^{\frac{d}{2}} \norm{\sum_{\und{i}\in \n^d_m}x_{\und{i}}\otimes a_{p,i_1} \otimes \cdots
\otimes a_{p,i_d}}_{L_p(\M \bar{\otimes} \Gamma^{\otimes d}_{\pm1})}.
\end{align*}
We get the desired conclusion by Corollary \ref{lem-tensor} again.

\end{proof}

\begin{rem}{\rm
The conditions \eqref{symmetry-degree-2} and \eqref{symmetry-degree-d} imply $x_{ii} = 0$ for any $i$ and
$x_{\und{i}} = 0$ for any $\und{i} \in \n^d_m - \I^d_m$.
}
\end{rem}

\bibliographystyle{amsplain}
\providecommand{\bysame}{\leavevmode\hbox
to3em{\hrulefill}\thinspace}

\end{document}